\documentclass[12pt]{article}

\usepackage{amsmath}
\usepackage{amssymb}
\usepackage{amscd}
\usepackage[all]{xy}

\usepackage{amstext}
\usepackage{amsthm}
\usepackage{a4}
\usepackage[mathscr]{eucal} 
\usepackage{mathrsfs}

\numberwithin{equation}{section}

\newtheorem{theorem}{Theorem}[section]
\newtheorem{definition}[theorem]{Definition}
\newtheorem{proposition}[theorem]{Proposition}
\newtheorem{corollary}[theorem]{Corollary}
\newtheorem{lemma}[theorem]{Lemma}
\newtheorem{remark}[theorem]{Remark}
\newtheorem{remarks}[theorem]{Remarks}
\newtheorem{example}[theorem]{Example}

\newcommand{\cali}[1]{\mathscr{#1}} 
\newcommand{\Rc}{{\cali R}}
\newcommand{\Cc}{{\cali C}}
\newcommand{\Uc}{{\cali U}}
\newcommand{\Kc}{{\cali K}}
\newcommand{\Ec}{{\cali E}}
\newcommand{\Sc}{{\cali S}}

\newcommand{\lov}{{\rm lov}}
\newcommand{\vol}{{\rm volume}}
\renewcommand{\L}{{\cali L}}
\newcommand{\Lone}{{{\cali L}^1}}
\newcommand{\Loneloc}{{\cali L}^1_{loc}}
\newcommand{\Linfty}{{{\cali L}^\infty}}
\newcommand{\C}{\mathbb{C}}

\newcommand{\Z}{\mathbb{Z}}

\newcommand{\R}{\mathbb{R}}

\newcommand{\ddc}{{{\rm dd}^{\rm c}}}
\renewcommand{\d}{{\rm d}}
\newcommand{\supp}{{\rm supp}}
\newcommand{\pr}{{\rm pr}}

\title{Geometry of currents, intersection theory and dynamics of horizontal-like maps}
\author{Tien-Cuong Dinh and Nessim Sibony}

\begin{document}

\maketitle

\begin{abstract} 
We introduce a geometry on the cone of positive closed currents of bidegree $(p,p)$ 
and apply  it to define the intersection of such currents.
We also construct and study the Green currents and the equilibrium measure for
horizontal-like mappings. The Green currents satisfy some extremality properties. 
The equilibrium measure is invariant, mixing and has maximal entropy. 
It is equal to the intersection of the Green currents associated to the horizontal-like map and
to its inverse.
\end{abstract}

\section{Introduction}
\label{introduction}

In this paper we develop the theory of positive closed currents of any degree in
order to continue our exploration of dynamical systems in several variables, 
with emphasis
on systems not defined by rational maps. 

In \cite{DinhSibony1}, we developed the theory of polynomial-like
maps in higher dimension. Recall that a polynomial-like map is a holomorphic map
$f:U\rightarrow V$, with $U\Subset V\Subset \C^k$, and that $f$ is proper of topological degree $d_t>1$.
In some sense, such a map is expanding, but it has critical points in general.

Here, we consider horizontal-like maps in any dimension. Basically, a horizontal-like map is a holomorphic map
defined on a domain in $\C^k$, which is ``expanding'' in $p$ directions and 
``contracting'' in $k-p$ directions.
The expansion and contraction are of global nature, but the map is, in general, not uniformly 
hyperbolic in the dynamical sense
\cite{KatokHasselblatt}.
The precise definition is given in Section \ref{horiz_map}.

This situation has been already studied by Dujardin for $k=2$ with emphasis on biholomorphic maps \cite{Dujardin}. 
The study was developed in dimension 2 by Dujardin and the authors to deal with the random iteration of
meromorphic horizontal-like maps,
in order to study rates of escape to infinity for polynomial mappings 
in $\C^2$ \cite{DinhDujardinSibony}. It turns out that, as for polynomial-like maps,  the building blocs 
for a large class of
polynomial maps are horizontal-like maps. We should observe that to treat the case 
of $\C^2$ with the methods of the present paper one should deal with horizontal-like
maps in $\C^4$ or $\C^8$ and that we obtain new results even in the $\C^2$ case
(see Theorem \ref{non_closed_current} and \cite{Dinh}). The main technical problem is to deal with currents of
higher bidegree.

One of the difficulties is that the potentials of currents of higher bidegree are not functions.
Hence, the techniques used in the case of dimension 2 do not work for general horizontal-like maps.
It seems that considering the potentials is not the best way to prove properties of currents of higher bidegree. 
We propose here
another approach to deal directly with the cone of positive closed currents that we consider as a space
of infinite dimension with some plurisubharmonic (p.s.h.) structure. 

We introduce in Section \ref{geometry_currents}
the notion of {\it structural varieties} in this cone which allows us to use the complex
structure of $\C^k$. 
Structural varieties connect currents in this cone. So, we will 
study singular currents using their smooth approximation in structural
discs. 
For example, a structural disc of currents of bidimension $(p,p)$ 
is the collection of slices 
of a positive closed 
current of bidimension $(p+1,p+1)$. The family is not always continuous in term
of slices, but when it acts on forms $\Phi$ such that $\ddc\Phi\geq 0$ we get p.s.h. 
functions on the space parametrizing slices. To prove the 
convergence of a sequence of currents we embed it in some sequence
of structural discs passing through a common smooth current. 
We then use systematically 
the convergence properties of the sequence of p.s.h. functions produced by the action 
on a test form $\Phi$ with $\ddc\Phi\geq 0$. 
An analog $\Lambda_\Phi$ of the Abel-Radon transform is also introduced. It 
plays the role of {\it p.s.h. functions} on the space of currents. 

In Section \ref{intersection_currents}, 
we use the structural discs in order to define the wedge product $T\wedge S$ where $T$ is a vertical
positive closed current and $S$ is a horizontal one, of the right bidegrees, 
such that the supports intersect on a compact set. 
Let $\varphi$ be a p.s.h. function on a small neighbourhood $W$ of $\supp(T)\cap\supp(S)$. We define
\begin{eqnarray}
\langle T\wedge S, \varphi \rangle &:=& \limsup_{T'\rightarrow T\atop S'\rightarrow S} 
\langle T'\wedge S',\varphi\rangle \label{eq_intersection}
\end{eqnarray}
where $T'$ are smooth vertical currents approaching $T$ and $S'$ are smooth 
horizontal currents approaching $S$ with $\supp(T')\cap\supp(S')\subset W$. We use structural discs in order to
show that the right hand side of (\ref{eq_intersection}) depends linearly on $T$, $S$ and $\varphi$.
This wedge product has
interesting continuity properties. 

We believe that the notion of structural discs
will be useful in other situations. It is a notion of deformation of a 
positive closed current into another one in the same ``homology'' class. 
This can be also useful in the context of compact manifolds. 

We apply the above theory of currents to study horizontal-like maps in
$\C^k$, $k\geq 2$. 
A horizontal-like map has a (main) dynamical degree $d$; this allows
us  
to define an operator $\L_v:=\frac{1}{d}f^*$ (resp. $\L_h:=\frac{1}{d}f_*$) 
on vertical (resp. horizontal) currents. One of our main results
is the following (Theorem \ref{random_green_current}). Let $f_n$ be a sequence of invertible horizontal-like maps
and let $R_n$ be a sequence of normalized vertical positive closed forms. If $R_n$ are uniformly bounded,
then 
$\L_{v,1}\ldots \L_{v,n}(R_n)$ converge to a normalized vertical current $T_+$ 
which is independent of $(R_n)$. 
If the $R_n$'s are continuous, the convergence is stronger than weak convergence (see Remark \ref{rk_random_test}).
We use structrural discs in the proof in order to deduce the convergence of currents from
the convergence of subharmonic functions on structural discs.

When all the $f_n$'s are equal to $f$, 
we obtain a Green current satisfying $f^*(T_+)=dT_+$ 
(Corollary \ref{green_current}). We are then able to produce in this case a mixing invariant measure $\mu$
(Theorem \ref{equilibrium_measure}).
This is done by going to the product space and applying our formalism to the horizontal-like map
$F:(x_1,x_2)\mapsto(f(x_1), f^{-1}(x_2))$ of dynamical degree $d^2$. 
More precisely, if $R$ (resp. $S$) is a normalized smooth
vertical (resp. horizontal) positive closed form then the equilibrium measure is constructed
as $\mu:=\lim d^{-2n}(f^n)^*R\wedge (f^n)_*S$. Formally, if $\Delta$ is the diagonal of the product space, we 
obtain $\mu$ as the limit of 
$$d^{-2n}\big((f^n)^*R\otimes (f^n)_*S\big)\wedge [\Delta]=d^{-2n} F^{n*}(R\otimes S)\wedge [\Delta].$$
This reduces the problem to the study of 
strong convergence of the vertical currents $d^{-2n}F^{n*}(R\otimes S)$
(see Remark \ref{rk_random_test}).  
We finally show that $\mu=T_+\wedge T_-$ (Theorem \ref{intersection_green_currents}).
Here, $T_-$ is the Green current associated to $f^{-1}$.

Our proof of the mixing of the equilibrium measure uses also a new 
idea different from the approach in 
Bedford-Smillie \cite{BedfordSmillie} for H{\'e}non maps or 
in \cite{Sibony} for regular polynomial automorphisms. The method is to use 
the maps of type $(x_1,x_2)\mapsto(f(x_1), f^{-1}(x_2))$ in order to reduce the
problem to a linear one.

Using classical arguments
\cite{Gromov, Yomdin, Walters, BedfordSmillie, DinhSibony1}, we show that $\mu$ has maximal entropy $\log d$
(Theorem \ref{entropy}).

\section{Geometry of currents}
\label{geometry_currents}

In this Section we study the geometry of the cones of positive closed currents which are 
supported in vertical or horizontal subsets of a domain $D=M\times N$. We 
define structural discs, p.s.h. functions and the Kobayashi pseudo-distance on these cones.
We refer to \cite{Federer, Lelong, Demailly, HarveyShiffman}
for the basics on the theory of currents.  
For the reader's convenience,
we recall some properties, that we use in this article,   
of the slicing operation in the complex
setting. 
\bigskip

\noindent
{\bf $\bullet$ Slicing theory.} 
Let $X$, $V$ be two complex manifolds of dimension $k+l$ and $l$ respectively. Let 
$\Pi_V:X\rightarrow V$ be a surjective holomorphic map and $\Rc$ be a current on $X$ of degree $2k+2l-m$ 
and of dimension $m$ with
$m\geq 2l$. Assume that  $\Rc$, $\partial \Rc$ 
and $\overline\partial \Rc$ are of order $0$. 
One can define {\it the slice} $\langle \Rc,\Pi_V,\theta\rangle$
 for almost every $\theta\in V$.
This is a current of dimension $m-2l$ on $\Pi_V^{-1}(\theta)$.
One can of course consider it as a current on $X$.
When $\Rc$ is of bidimension $(n,n)$, $\langle \Rc,\Pi_V,\theta\rangle$ are
of bidimension $(n-l,n-l)$. The slicing commutes with the operations $\partial$
and $\overline\partial$. In particular, if $\Rc$ is closed then $\langle \Rc,\Pi_V,\theta
\rangle$ is also closed. 

Slicing is the generalization of restriction of forms to level sets of
holomorphic maps.
When $\Rc$ is a continuous form, 
$\langle \Rc, \Pi_V,\theta\rangle$ is simply the restriction of $\Rc$ to $\Pi_V^{-1}(\theta)$. 
When $\Rc$ is the current of integration on an analytic subset $Y$ of $X$, 
$\langle \Rc,\Pi_V,\theta\rangle$ is the current of integration on the analytic set 
$Y\cap \Pi_V^{-1}(\theta)$ for $\theta$ generic. If $\varphi$ is a continuous form on $X$ then
$\langle \Rc\wedge\varphi,\Pi_V,\theta\rangle=\langle \Rc,\Pi_V,\theta\rangle\wedge\varphi$.     

Let $y$ denote the coordinates in a chart of $V$ and $\lambda_V$ the standard volume form.
Let $\psi(y)$ be a positive smooth function with compact support
such that $\int\psi\lambda_V=1$. Define
$\psi_\epsilon(y):=\epsilon^{-2l}\psi(\epsilon^{-1} y)$ and 
$\psi_{\theta,\epsilon}(y):=\psi_\epsilon(y-\theta)$ (the measures $\psi_{\theta,\epsilon}\lambda_V$
approximate the Dirac mass at $\theta$). Then, for every smooth test form $\Psi$
of the right degree with compact support in $X$ 
one has
$$\langle \Rc,\Pi_V,\theta\rangle (\Psi)=\lim_{\epsilon\rightarrow 0}
\langle \Rc\wedge \Pi_V^*(\psi_{\theta,\epsilon}\lambda_V),\Psi\rangle$$
when $\langle \Rc,\Pi_V,\theta\rangle$ exists.
This property holds for all choice of the function $\psi$ and for $\Psi$ such that $\Pi_V$
is proper on $\supp(\Psi)\cap\supp(\Rc)$.
Conversely, when the previous limit exists and is independent of $\psi$, 
it defines $\langle \Rc,\Pi_V,\theta\rangle$ and one says 
that $\langle \Rc,\Pi_V,\theta\rangle$ {\it is well defined}.
We have the following formula for every continuous form $\Omega$ 
of maximal degree with compact support in $V$:
\begin{eqnarray} \label{eq_property_slice}
\int_V\langle \Rc,\Pi_V,\theta\rangle (\Psi)\Omega(\theta) & = &  \langle \Rc\wedge \Pi_V^*(\Omega),
\Psi\rangle.
\end{eqnarray}

We will show that in the situation we consider, slices are always well defined.

\bigskip

\noindent
{\bf $\bullet$ Vertical and horizontal currents.} 
Let $M\subset \C^p$ and $N\subset \C^{k-p}$ be two bounded 
convex open sets (see Remark \ref{rk_nonconvex}). Consider the 
domain $D:=M\times N$ in $\C^k$. We call {\it vertical} 
(resp. {\it horizontal}) 
{\it boundary} of $D$ the set $\partial_v D:=\partial M\times N$ 
(resp. $\partial_h D:=M\times \partial N$). 
A subset $E$ of $D$ is called {\it vertical} (resp. {\it horizontal}) if
$\overline E$ does not intersect  $\overline{\partial_v D}$ 
(resp. $\overline{\partial_h D}$). Let $\pi_1$ and $\pi_2$ denote the canonical 
projections of $D$ on $M$ and $N$. Then, $E$ is vertical (resp. horizontal) if and only if 
$\pi_1(E)\Subset M$ (resp. $\pi_2(E)\Subset N$).
A current on $D$ is {\it vertical} (resp. 
{\it horizontal}) if its support is vertical (resp. horizontal). 

Let $\Cc_v(D)$ (resp. $\Cc_h(D)$) denote 
the cone of {\bf positive closed} vertical
(resp. horizontal) currents of bidegree $(p,p)$ (resp. $(k-p,k-p)$) on $D$.   
Consider a current $R$ in $\Cc_v(D)$. 
Since $\pi_2$ is proper on $\supp(R)$, $(\pi_2)_*(R)$ is a positive closed current of
bidegree $(0,0)$ on $N$. Hence, $(\pi_2)_*(R)$ is given by a constant function $c$ on
$N$. 
Formula (\ref{eq_property_slice}) implies that 
the mass of the slice measure $\langle R,\pi_2,w\rangle$
is independent of $w$ and is equal to $c$. 
We will show in Theorem \ref{th_slices} that in this 
situation, the slice measure is defined for every $w\in N$
(see also Theorem \ref{th_intersection_currents}). 
We say that $c$ is 
the {\it slice mass} of $R$ and we denote it by $\|R\|_v$.
For every smooth probability measure $\Omega$ with compact support in $N$, we have 
$\|R\|_v:=\langle R,(\pi_2)^*(\Omega)\rangle$.  
When $\|R\|_v=1$ we say that $R$ is {\it normalized}. Let $\Cc_v^1(D)$ denote the set 
of such currents.
The slice mass $\|\cdot\|_h$ and the convex $\Cc_h^1(D)$ for 
horizontal currents are similarly defined.

\bigskip

\noindent
{\bf $\bullet$ Structural varieties and p.s.h. functions.} 
In order to use the complex structure of $D$, we introduce the notion 
of structural varieties in $\Cc^1_v(D)$. Let $V$ be a connected complex manifold.
Let $\Rc$ be a positive closed current of bidegree $(p,p)$ in $V\times D$. 
Let $\Pi_V:V\times D\rightarrow V$,  $\Pi_D:V\times D\rightarrow D$, $\Pi_M:V\times D\rightarrow M$ and
$\Pi_N:V\times D\rightarrow N$  
be the canonical projections. 
We assume that for every compact set $K\subset V$ the projection of 
$\supp(\Rc)\cap \Pi_V^{-1}(K)$ on $M$ is relatively compact in $M$. In particular
$\supp(\Rc)\cap \Pi_V^{-1}(\theta)$ is a vertical set of $\{\theta\}\times D$ 
for every $\theta\in V$. 

\begin{theorem} \label{th_slices}
For every $\theta\in V$ the slice $\langle \Rc, \Pi_V, \theta \rangle$ exists and is a vertical positive
closed current on $\{\theta\}\times D$. 
Moreover its slice mass is independent of $\theta$. If $\Psi$ is a real continuous 
$(k-p,k-p)$-form on $V\times D$ such that
 $\ddc\Psi\geq 0$  and  $\Pi_N(\supp(\Psi))\Subset N$ 
then
$\langle \Rc, \Pi_V, \theta\rangle (\Psi)$ defines a p.s.h. function on $V$. If $\ddc\Psi=0$ then
   $\langle \Rc, \Pi_V, \theta\rangle (\Psi)$ is pluriharmonic.
\end{theorem}
\begin{proof}
The problem is local, so we can assume that $V$ is a ball.
Consider the current $\Rc':=\Rc\wedge \Psi$ of bidegree $(k,k)$  on $V\times D$. 
It satisfies $\ddc \Rc'\geq 0$.
Observe that for every 
$\theta\in V$, 
$\supp(\Rc')\cap \Pi_V^{-1}(\theta)$ is compact in $\{\theta\}\times D$
and $\Pi_V$ is proper on the support of $\Rc'$. Then 
$(\Pi_V)_*(\Rc')$ is well defined. It is a current of bidegree $(0,0)$
on $V$ which satisfies $\ddc (\Pi_V)_*(\Rc')\geq 0$. 
Therefore, it is defined 
by a p.s.h. function $\varphi$.
It follows that if $\psi$, $\psi_{\theta,\epsilon}$ and $\lambda_V$ are as above then 
$\int \varphi \psi_{\theta,\epsilon}\lambda_V$ converges to $\varphi(\theta)$. 

The last integral is equal to $\langle \Rc\wedge \Pi_V^*(\psi_{\theta,\epsilon}\lambda_V),\Psi\rangle$.
Hence 
$\langle \Rc\wedge
\Pi_V^*(\psi_{\theta,\epsilon}\lambda_V),\Psi\rangle$ converge to
$\varphi(\theta)$ which does
not depend on $\psi$. 
This holds also for every  smooth test form $\Psi'$ with compact
support in $V\times D$. Indeed, we have the following lemma.

\begin{lemma} \label{lemma_density}
Smooth $(k-p,k-p)$-forms with compact support in $V\times D$ 
belong to the space generated by the forms $\Psi$ satisfying the
hypotheses of Theorem \ref{th_slices}, i.e. $\Pi_N(\supp(\Psi))\Subset
N$ and $\ddc\Psi\geq 0$.
\end{lemma}
\proof
Let $\Psi'$ be a real smooth $(k-p,k-p)$-form with compact support in $V\times D$. Let
$\Omega$ be a positive form of maximal degree on $N$, with compact support
and strictly positive on $\Pi_N(\supp(\Psi'))$. If $\rho$ is
a smooth strictly p.s.h. function on $V\times D$ then $\Psi_0:=\rho\Pi_D^*(\pi_2^*(\Omega))$
is a smooth form satisfying the hypothesis of Theorem \ref{th_slices}. If
$\pi_{2,\epsilon}$ is a small pertubation of $\pi_2$ then
$\Psi_\epsilon:=\rho\Pi_D^*(\pi_{2,\epsilon}^*(\Omega))$ 
satisfies the same properties. Taking a linear combination of such
forms we obtain a form $\Psi$ such that $\ddc\Psi$ is strictly
positive on $\supp(\Psi')$. Then we can write
$\Psi'=(A\Psi+\Psi')-A\Psi$ with $A>0$ large enough. The forms
$A\Psi+\Psi'$ and $A\Psi$ satisfy  the
hypotheses of Theorem \ref{th_slices}, in particular, we have 
$\ddc(A\Psi+\Psi')\geq 0$ and $\ddc(A\Psi)\geq 0$. 
\endproof

Hence $\langle \Rc, \Pi_V,\theta\rangle$ is well defined and 
 $\langle \Rc, \Pi_V,\theta\rangle(\Psi)=\varphi(\theta)$ is a p.s.h. function on $\theta$.
When $\ddc\Psi=0$ the function $\langle \Rc, \Pi_V,\theta\rangle(-\Psi)$ is also p.s.h.
Hence $\langle \Rc, \Pi_V,\theta\rangle(\Psi)$ is pluriharmonic.

Let $\Omega$ be as above. Consider $\Psi:=\Pi_D^*(\pi_2^*(\Omega))$. In this case, since $\Psi$ is closed,
$\varphi$ is also closed. It follows that $\varphi$ is a constant function. By definition 
$\varphi(\theta)=\langle \Rc,\Pi_V,\theta\rangle(\Psi)$ is equal to the slice mass of 
$\langle \Rc,\Pi_V,\theta\rangle$. Therefore, the slice mass is independent of
$\theta$.
\end{proof}

\begin{remark} \label{rk_slices} \rm
One can identify $R_\theta=\langle \Rc,\Pi_V,\theta\rangle$ with a current in $\Cc_v(D)$. 
Theorem \ref{th_slices} implies that the family $(R_\theta)$ is
continuous for the {\it plurifine topology} on $V$,
i.e. the coarsest topology for which p.s.h. functions are continuous. 
Let $\Phi$ be a real horizontal current of bidegree $(k-p,k-p)$, of finite 
mass on $D$ such that $\ddc\Phi\geq 0$. 
If $\Rc$ or $\Phi$ is a continuous form then $\langle \Rc,\Pi_V, \theta\rangle (\Phi)$
defines a p.s.h. function on $V$. Indeed, we can apply Theorem \ref{th_slices} to $\Psi:=\Pi_D^*(\Phi)$. 
\end{remark}

\begin{definition} \label{def_structural_variety}
\rm
Theorem \ref{th_slices} allows us to define a map
$\tau:V\rightarrow \Cc_v(D)$
$$\tau(\theta):=R_\theta=\langle \Rc,\Pi_V,\theta\rangle.$$ 
If we multiply $\Rc$ by a suitable constant, all the values of $\tau$ are normalized.
We say that $\tau$
defines a {\it structural variety} in $\Cc^1_v(D)$. 

A function $\Lambda:\Cc_v^1(D)\rightarrow 
\R\cup\{-\infty\}$ is called {\it p.s.h.} if it is not identically equal to $-\infty$ and for every 
structural variety $\tau:V\rightarrow \Cc_v^1(D)$ the function $\Lambda\circ\tau$ is either p.s.h. or 
identically $-\infty$ on $V$. If $\Lambda$ and $-\Lambda$ are p.s.h. we say that $\Lambda$
is {\it pluriharmonic}.
\end{definition}

Let $\Phi$ be a real continuous horizontal  $(k-p,k-p)$-form on $D$. 
Define the linear map $\Lambda_\Phi:\Cc_v(D)\rightarrow \R$
by $\Lambda_\Phi(R):=\langle R,\Phi\rangle$.
Such an operator is a version of the Abel-Radon
transform in complex analysis.
Observe that
real smooth $(k-p,k-p)$-forms with compact support in $D$ belong to the space generated 
by the smooth horizontal forms $\Phi\geq 0$ with $\ddc \Phi\geq 0$
(see Lemma \ref{lemma_density}).
Hence, the maps $\Lambda_\Phi$ with $\Phi\geq 0$ and $\ddc\Phi\geq 0$, 
separate currents in $\Cc_v(D)$. 
Theorem \ref{th_slices} and Remark \ref{rk_slices} show
that $\Lambda_\Phi$ is p.s.h. on $\Cc_v^1(D)$ when $\ddc\Phi\geq 0$.

We can summarize our construction of the function $\theta\mapsto \langle \tau(\theta),\Phi\rangle$
by the following diagram:
\begin{eqnarray}
\xymatrix{
          V \ar[r]^-{\tau}  & \Cc^1_v(D)    \ar[r]^-{\Lambda_\Phi} & \R\\
   } \label{diagram1}
\end{eqnarray}

\noindent
$\bullet$ {\bf Some structural discs.} Given a current $R\in \Cc_v^1(D)$, we construct some special 
structural discs passing through $R$, that we will use
in the next sections. For these discs, the map $\tau$ 
is continuous with respect to the weak topology on currents. 
In order to construct the current $\Rc$,
we consider the images of 
$R$ under holomorphic families of maps.

Let $M'\Subset M$ and $N''\Subset N$ be open sets. Define 
$D':=M'\times N$ and $D'':=M\times N''$. 
In order to simplify the notations, assume
that $0$ belongs to $(M\setminus\overline M')\times (N\setminus\overline N'')$. 
Fix a domain $D^*=M^*\times N^*\Subset  D$ with $M\setminus M^*$ and
$N\setminus N^*$ small enough, $M'\Subset M^*$, $N'\Subset N^*$.
Choose  
a small simply connected neighbourhood $V$ of $[0,1]$ in $\C$.
Finally, choose a small open neighbourhood
$U\Subset D^*$ of $0$ in $\C^p\times \C^{k-p}$ and a smooth positive function $\rho$ with
support in $U$ such that $\int \rho(a,b)\lambda(a,b)=1$. Here, $\lambda$
denotes the standard volume form on $\C^k$.

For $\theta\in V$ and $(a,b)\in U$,  
define the affine map
$h_{a,b,\theta}:\C^p\times\C^{k-p}\rightarrow \C^p\times\C^{k-p}$ by
$$h_{a,b,\theta}(z,w):= \big(\theta z+(1-\theta)a,
w+(\theta-1)b\big).$$
These maps are small pertubations of the map $(z,w)\mapsto (\theta z,w)$.
When $\theta=1$ we obtain the identity map and when $\theta=0$ we obtain an affine
map onto the subspace $\{z=a\}$. 
Let $R$ be a 
current in $\Cc^1_v(D')$. We will show that the currents 
$R_{a,b,\theta}:= (h_{a,b,\theta})_*(R)$ define a structural disc in 
$\Cc^1_v(D^*)$, i.e. they are slices of a current $\Rc_{a,b}$ in $V\times D^*$. 

Observe that $R_{a,b,\theta}$ is well defined, since 
$h_{a,b,\theta}:\supp(R)\cap h_{a,b,\theta}^{-1}(D^*)\rightarrow
D^*$ is proper. This last  property follows from
the fact that $M$ is convex and  $h_{a,b,\theta}$
is close to the map $(z,w)\mapsto (\theta z, w)$.
Moreover, $R_{a,b,\theta}$ is well defined on some open set $D_\theta$
which converges to $D$ when $\theta\rightarrow 1$.
The dependence of currents $R_{0,0,\theta}$ on $\theta$ has 
been used by Dujardin in order to study H{\'e}non-like maps \cite{Dujardin}
(see also \cite{DinhDujardinSibony}).

Define the meromorphic map
$H_{a,b}:V\times D^*\rightarrow \C^p\times N$ 
by 
$$H_{a,b}(\theta,z,w):=h_{a,b,\theta}^{-1}(z,w)=\left(\frac{z+ 
(\theta-1)a}{\theta}, w-(\theta-1)b\right).$$ 
The current
$\Rc_{a,b}:=H_{a,b}^*(R)$, which is of bidimension $(k-p+1,k-p+1)$, is well defined out of the pole set 
$\{\theta=0\}$ of $H_{a,b}$. Since $\supp(\Rc_{a,b})\subset H_{a,b}^{-1}(\supp(R))$, then
when $\theta$ approaches 0, $\supp(\Rc_{a,b})$ clusters
only on the set $\{z=a\}$. So,
this current is well defined out of $\{\theta=0\}\cap\{z=a\}$.
The dimension of $\{\theta=0\}\cap\{z=a\}$, which is equal to $k-p$, is smaller than 
the dimension of $\Rc_{a,b}$. Hence,
one can extend $\Rc_{a,b}$  across  
$\{\theta=0\}\cap \{z=a\}$ with no mass on this set \cite{HarveyPolking}. 

Since $M$ is convex and since
$h_{a,b,\theta}$ is close to the map $(z,w)\mapsto (\theta z,w)$,
$\supp(\Rc_{a,b})\cap\Pi_V^{-1}(\theta)$, which is isomorphic to 
$\supp(R_{a,b,\theta})$, 
is a vertical set of $\{\theta\}\times D^*$ 
for every $\theta\in V$. 
Hence, the slice currents $\langle \Rc_{a,b},\Pi_V, \theta\rangle$ 
define a structural 
disc in $\Cc^1_v(D^*)$. By Theorem \ref{th_slices},
these slices exist for every $\theta\in V$
and are equal to $R_{a,b,\theta}$ (we identify $\{\theta\}\times D$ with $D$). 
The currents $R_{a,b,\theta}$ depend continuously on $\theta$ for the 
weak topology on currents. This is clear for $\theta\not=0$, and as we have seen, 
the limit at $\theta=0$  is $[z=a]$
(see also Lemmas \ref{disc1} and \ref{disc2} below).
Recall that $[z=a]$ denotes the current of integration on 
the analytic set $\{z=a\}$.

We have
$R_{a,b,1}=R$ and 
$R_{a,b,0}=[z=a]$.
Hence, $R_{a,b,0}$ is independent of $R$.
In other words, when $R$ varies we obtain a family of structural discs passing through the 
same point 
$[z=a]$ in $\Cc^1_v(D^*)$.

We introduce a smoothing. Define 
$$\Rc:=\int \Rc_{a,b}\rho(a,b) \lambda(a,b).$$
The current $\Rc$ satisfies the hypothesis of Theorem \ref{th_slices} for $D^*$.
Hence, the slice currents $R_\theta:=\langle \Rc,\Pi,\theta\rangle$
define a structural disc in $\Cc^1_v(D^*)$. These slices are well
defined for every $\theta\in V$ and 
\begin{eqnarray}
R_\theta=\int R_{a,b,\theta}\rho(a,b) \lambda(a,b). \label{eq_R_theta}
\end{eqnarray}
Observe that $R_\theta$ depends continuously on $\theta$ for the weak topology. We have 
$R_1=R$ and 
$$R_0=\int [z=a]\rho(a,b)\lambda(a,b)=\pi_1^*(\pi_1)_*(\rho\lambda).$$ 
The last current 
is independent of $R$. When $R$ varies, we obtain again a family of structural discs which pass
through the same point $\pi_1^*(\pi_1)_*(\rho\lambda)$ in $\Cc^1_v(D^*)$.

In the following two lemmas, we study the continuity of $R_\theta$ near $0$ and near $1$. We will use 
these facts in our 
convergence theorems. Lemma \ref{disc1} shows that every current in $\Cc^1_v(D')$ can be joined to a fixed 
vertical current $R_0$ through smooth ones.

\begin{lemma} \label{disc1} Let $R\in\Cc^1_v(D')$. 
Then, for $\theta\in V\setminus\{1\}$, 
$R_\theta$ is a smooth form on $D^*$. For
  $\theta\in V\setminus\{0,1\}$, 
$R_\theta$ depends continuously on $(R,\theta)$
for the $\Cc^\infty$ topology. 
Moreover,
there exist $r>0$ and $c>0$ independent
of $R$ such that if $|\theta|\leq r$
$$\|R_\theta-R_0\|_{\Linfty(D^*)}\leq c|\theta|$$
where the $\Linfty$ norm on forms is the sum of $\Linfty$ norms of coefficients.
\end{lemma}
\begin{proof} 
The smoothness of $R_\theta$ for 
$\theta\not =1$, and their dependence of $(R,\theta)$
are checked using a classical change of variables in (\ref{eq_R_theta}) as follows.
Let $\widetilde R$, $\widetilde R_{a,b,\theta}$ and $\widetilde R_\theta$ be 
the coefficients of $\d z_I\wedge \d \overline z_J\wedge \d w_K\wedge \d \overline w_L$ in 
$R$, $R_{a,b,\theta}$ and $R_\theta$ respectively, for some multi-indices $I$, $J$, $K$, $L$. Let 
$$(A,B):= h^{-1}_{a,b,\theta}(z,w) = \left(\frac{z+ 
(\theta-1)a}{\theta}, w-(\theta-1)b\right)$$
be the new variables. Since $R_{a,b,\theta}=(h_{a,b,\theta}^{-1})^* R$ we have
$$\widetilde R_{a,b,\theta}(z,w) =\theta^{-|I|} \overline
\theta^{-|J|} \widetilde R(A,B)$$ 
and from (\ref{eq_R_theta})
$$\widetilde R_\theta(z,w) =\theta^{-|I|} \overline
\theta^{-|J|} \int_{A,B} \widetilde R(A,B)\cdot (\rho\lambda)\left( \frac{\theta A -z}{\theta-1},
\frac{w-B}{\theta-1}\right).$$ 
The smoothness of $R_\theta$ for $\theta\in V\setminus\{0,1\}$ and the dependence of $(R,\theta)$
are clear.

Let $\Phi$ be a smooth $(k-p,k-p)$-form  with compact support 
in $D^*$. By duality, the inequality that we have to prove is equivalent to
$$|\langle R_\theta -R_0,\Phi\rangle| \leq c|\theta|\|\Phi\|_\Lone.$$ 
From (\ref{eq_R_theta}), we get
\begin{eqnarray*}
\langle R_\theta,\Phi\rangle & = & 
\int \langle R_{a,b,\theta},\Phi \rangle \rho(a,b) \lambda(a,b)
= \int\langle R, h_{a,b,\theta}^*(\Phi) \rangle \rho(a,b) \lambda(a,b)\\
& = & \left\langle R, \int h_{a,b,\theta}^*(\Phi)\rho(a,b) \lambda(a,b) 
\right \rangle =: \langle R,\Phi_\theta\rangle.
\end{eqnarray*}
This also holds for $\theta=0$ by continuity. 
The forms $\Phi_\theta$ are obtained by convolution. They are smooth and uniformly bounded
by $c\|\Phi\|_\Lone$. Using the change of variables $(a,b)\mapsto (A,B):=h_{a,b,\theta}(z,w)$, we get
$$\|\Phi_\theta-\Phi_0\|_{\Linfty}\leq c |\theta| \|\Phi\|_\Lone$$
for $\theta$ small. Lemma \ref{disc1} follows.
\end{proof}

\begin{remark} \label{rk_nonconvex} \rm
When $M$ is not convex and $V'\subset V$ is a small neighbourhood of
1, then $(R_\theta)_{\theta\in V'}$ defines also a structural disc in
$\Cc^1_v(D'')$. The first part of Lemma \ref{disc1} holds in this case.
\end{remark}

\begin{lemma} \label{disc2} Let $R\in \Cc^1_v(D')$ be a continuous form.
Let $m(R,\epsilon)$ denote the modulus of continuity of $R$. Then, 
there exist $r>0$, $c>0$, $A>0$ independent of $R$ such that for 
$|\theta-1|\leq r$
$$\|R_\theta-R\|_{\Linfty(D^*)}\leq c \big(\|R\|_{\Linfty(D)}|\theta-1| +
m(R,A|\theta-1|)\big).$$
\end{lemma}

\begin{proof} Let $W$ denote the disc $\{|\theta-1|\leq r\}$ with $r>0$
small enough, so we are away of $\{\theta=0\}$. Then, there exists $A>0$ such that 
$\|h_{a,b,\theta}^{-1}(z,w)-(z,w)\|_{\Cc^1}\leq A |\theta-1|$ 
when $(z,w,a,b,\theta)\in D\times U\times W$.
Hence, there exists $c>0$ such that for $(a,b)\in U$
$$\|R_{a,b,\theta}-R\|_{\Linfty(D^*)}
\leq c\big(\|R\|_{\Linfty(D)} |\theta-1|+
m(R,A|\theta-1|)\big).$$
We can also prove this inequality using the description of $\widetilde R_{a,b,\theta}$ as 
in Lemma \ref{disc1}.  
Finally, we obtain the desired inequality by integration using (\ref{eq_R_theta}).
\end{proof}

\noindent
{\bf $\bullet$ Kobayashi pseudo-distance.} 
Let $\Cc_v^1(\overline D)$ be the set of currents in $\Cc_v^1(D)$ which can be extended 
to a current in $\Cc_v^1(M\times N')$ for some neighbourhood $N'$ of $\overline N$. 
We introduce {\it the Kobayashi pseudo-distance} $\rho_v$ on $\Cc_v^1(\overline D)$. 
Let $R$ and $S$ be two currents in $\Cc_v^1(\overline D)$.
Let $\Delta$ be the unit disc and $\rho_0$ denote the hyperbolic distance on $\Delta$.
Consider chains of continuous 
structural discs $\tau_i:\Delta\rightarrow \Cc_v^1(\overline D)$ which connect $R$ and $S$.
More precisely, suppose $\theta_i$, $\theta_i'\in \Delta$ such that
$\tau_1(\theta_1)=R$, $\tau_i(\theta_i')=\tau_{i+1}(\theta_{i+1})$ and $\tau_n(\theta_n')=S$. 
Define
\begin{eqnarray}
\rho_v(R,S) &:= & \inf \sum_{i=1}^n \rho_0(\theta_i,\theta_i') \label{def_pseudodistance}
\end{eqnarray} 
where the infimum is taken over all the $n$, $\tau_i$, $\theta_i$ and $\theta_i'$. We have seen that $R$ and $S$ can
be connected by a chain of 
two continuous structural discs. Hence $\rho_v(R,S)$ is finite. It is easy to check that $\rho_v$
satisfies the triangle inequality.

\begin{proposition} \label{pr_hyperbolicity}
The pseudo-distance $\rho_v$ is not a distance. If real continuous horizontal $\ddc$-closed forms on $D$
separate $R$ and $S$ then 
$\rho_v(R,S)>0$.
\end{proposition}
\begin{proof}
We have to construct two different currents $R$ and $S$ such that $\rho_v(R,S)=0$. 
We can replace $N$ by a ball containing $N$ and $M$ by a polydisc contained in $M$.
So, we can assume that $N$ is the unit ball and $M$ is the unit polydisc.
It is sufficient to consider the case where $p=1$ and $M$ is the unit disc $\Delta$. We obtain the general case 
by taking the product of  $M$ and $D$ by $\Delta^{p-1}$. 

Let $\nu_r$ be the Lebesgue measure on the circle $\{|z|=r\}$
normalized by $\|\nu_r\|=1$. 
Consider $R:=\pi_1^*(\nu_0)=[z=0]$ and $S:=\pi_1^*(\nu_{1/2})$. Let $\Rc$ be the positive closed current
of bidegree $(1,1)$ on $\Delta\times (\Delta\times \C^{k-1})$ 
given by $\Rc:=\ddc \Uc$ where $\Uc(\theta,z,w):=
\max\{\log|z|, \frac{1}{A}\log|\theta|\}$
and $A>1$. This current has support in $\{|z|^A=|\theta|\}$. Hence, if $\Pi_\Delta$ is the projection
on the first factor $\Delta$, 
the slices $R_\theta:=\langle\Rc,\Pi_\Delta,\theta\rangle=\ddc \Uc(\theta,\cdot,\cdot)$ 
define a continuous structural disc in $\Cc_v^1(\overline D)$.
Moreover, we have $R_\theta=\pi_1^*(\nu_r)$ where $r^A=|\theta|$.
In particular, we have $R_0=R$ and $R_\theta=S$ for $\theta:=2^{-A}$. It follows that 
$\rho_v(R,S)\leq \rho_0(0,\theta)$.  When $A\rightarrow\infty$, we have $\theta\rightarrow 0$ and then
$\rho_0(0,\theta)\rightarrow 0$. Therefore, $\rho_v(R,S)=0$. 

Now assume that $R$, $S$ satisfy the hypothesis of Proposition \ref{pr_hyperbolicity} and 
consider structural discs $\tau_i$ 
as above. 
Let $\Phi$ be a real continuous horizontal form such that $\ddc\Phi=0$ and $\langle R,\Phi\rangle
\not = \langle S,\Phi\rangle$.  
Using a regularization we can assume that $\Phi$
is smooth and is defined on a neighbourhood of $D$. Hence there exists a smooth current $\Phi'$ 
in $\Cc_h(D)$ such that $-\Phi'\leq \Phi\leq \Phi'$. 
Using coordinate changes on $\Delta$, one can also assume that 
$\theta_i=0$. If $\rho_0(0,\theta_i')>1$ then the right hand side of (\ref{def_pseudodistance}) is larger than 1.
We have only to consider the case where $\rho_0(0,\theta_i')\leq 1$ for every $i$. 

Define $\psi_i:=\Lambda_\Phi\circ\tau_i$. Theorem \ref{th_slices} 
implies that these functions are harmonic.
Lemma \ref{mass_of_intersection} below implies that they are
uniformly bounded by $\pm\|\Phi'\|_h$. 
Hence by Harnack's inequality 
$|\psi_i(\theta_i')-\psi_i(0)|\leq c \rho_0(0,\theta_i')$, $c>0$. On the other hand, we have 
$\psi_1(0)=\langle R,\Phi\rangle$, $\psi_i(\theta'_i)=\psi_{i+1}(0)$ and $\psi_n(\theta_n')=\langle S,\Phi\rangle$.
We then deduce that the right hand side of (\ref{def_pseudodistance}) is bounded from below by 
$c^{-1}|\psi_1(0)-\psi_n(\theta_n')|= c^{-1}|\langle R,\Phi\rangle - \langle S,\Phi\rangle|$. Hence $\rho_v(R,S)>0$.
\end{proof}

\begin{proposition} \label{pr_hyperbolicity_bis}
The space $\Cc^1_v(\overline D)$ is hyperbolic in 
the sense of Brody. More precisely, 
there exists no non-constant structural line
$\tau:\C\rightarrow\Cc^1_v(\overline D)$.
\end{proposition}
\begin{proof}
Consider a horizontal positive 
test forms $\Phi$ such that $\ddc\Phi\geq 0$
and assume $\Phi\leq \Phi'$ with $\Phi'$ a smooth form in $\Cc_h(D)$. Then, 
$\Lambda_\Phi\circ\tau$ is
constant since, by Theorem \ref{th_slices} and Lemma \ref{mass_of_intersection} below, 
it is a 
subharmonic function on $\C$, bounded from above 
by $\|\Phi'\|_h$. Proposition \ref{pr_hyperbolicity_bis} follows.
\end{proof}

\noindent
$\bullet$ {\bf Case of bidegree (1,1).} Assume that $p=1$. We will construct an example
of non-continuous structural discs in $\Cc_v^1(D)$. 
Let $\Rc$ be a positive closed current
of bidegree $(1,1)$ on $\Delta\times D$ satisfying the hypotheses of Theorem \ref{th_slices}.  
We can write $\Rc=\ddc \Uc$ where $\Uc$ is 
a p.s.h. function on $\Delta\times D$ which is pluriharmonic near $\Delta\times \partial_v D$. The slice 
$R_\theta:=\langle \Rc,\Pi_\Delta,\theta\rangle$ is equal to 
$\ddc U_\theta$ where $U_\theta:=\Uc_{|\{\theta\}\times D}$. The geometry of the support of $\Rc$
insures that $U_\theta$ is not identically equal to $-\infty$. Hence slice currents exist for every $\theta$.

Let $v$ be a bounded 
subharmonic function on $\Delta$. 
In order to simplify the notation assume that $0$ belongs to $D$.
Consider the case where 
$$\Uc(\theta,z,w)=\max \{v(\theta)-A,\log|z|\}.$$
The constant $A$ is chosen large enough so that $\Uc=\log|z|$ near $\Delta\times\partial_vD$. Then 
$\Rc$ vanishes near    $\Delta\times\partial_vD$. One easily check that $R_\theta=\pi_1^*(\nu_r)$
where $r:=\exp(v(\theta)-A)$. Hence $(R_\theta)$ is 
continuous with respect to $\theta$, if and only if $v$ is continuous.

The following proposition gives 
the converse of Theorem \ref{th_slices} in the bidegree $(1,1)$ case.

\begin{proposition} \label{pr_slices_(1,1)}
Let $R_\theta$ be a family of
currents of bidegree $(1,1)$ in  $\Cc^1_v(\{\theta\}\times D)$, $\theta\in\Delta$. Assume that the projection of 
$\cup\supp(R_\theta)$ on $M$ is relatively compact in $M$. 
Assume also that 
for every real continuous $(k-1,k-1)$-form $\Psi$ on $\Delta\times D$ such that
$\ddc\Psi\geq 0$ and $\Pi_N(\supp(\Psi))\Subset N$, the function
$\theta\mapsto \langle R_\theta,\Psi\rangle$ is subharmonic on $\Delta$. Then $\theta\mapsto
R_\theta$ defines a structural disc in $\Cc_v^1(D)$.  
\end{proposition}
\begin{proof} We want  to construct a potential $\Uc$ of a current $\Rc$ with given slices 
$R_\theta$. We will obtain $\Uc$ as a decreasing limit of some p.s.h. functions 
$\Uc_{\epsilon,\delta}$.

Let $\lambda$ denote the canonical volume form on $\C^{k-1}$ and 
$\psi$ be a positive radial function with compact support in $\C^{k-1}$ such that $\int\psi\lambda=1$. 
Define continuous functions $\psi_\epsilon(w):=\epsilon^{2-2k}\psi(\epsilon^{-1} w)$, $\epsilon>0$, 
and $\log_\delta|z|:=\max\{\log |z|, \log \delta\}$,
$\delta > 0$. Define also
$$\Phi^{\epsilon,\delta}_{z_0,w_0}(z,w):=\log_\delta|z-z_0| \psi_\epsilon(w-w_0)\lambda(w-w_0)$$
which is a regularization of the current $\log|z-z_0|\cdot[w=w_0]$, and 
$$\Uc_{\epsilon,\delta}(\theta_0,z_0,w_0):= \langle R_{\theta_0}, 
\Phi^{\epsilon,\delta}_{z_0,w_0}\rangle.$$
Here we identify $R_{\theta_0}$ to a current on $D$.

We first prove that for every domain $N^*\Subset N$, the function 
$\Uc_{\epsilon,\delta}$ is p.s.h. on $\Delta\times \C\times N^*$ for $\epsilon$ small enough. Assume that 
$z_0=g(\theta_0)$ and $w_0=h(\theta_0)$ where $(g,h)$ is a holomorphic map from $\Delta$ to $\C\times N^*$.
It is enough to prove
that $\Uc_{\epsilon,\delta}(\theta_0,g(\theta_0),h(\theta_0))$ is a subharmonic function with respect to $\theta_0$. 
This follows from the hypothesis. Indeed, in this case $\Phi_{z_0,w_0}^{\epsilon,\delta}(z,w)$ is equal to 
a continuous form $\Psi^{\epsilon,\delta}(\theta_0,z,w)$ which satisfies $\ddc\Psi^{\epsilon,\delta}\geq 0$
on $\Delta\times D$
and if $\epsilon$ is small enough $\Pi_N(\supp(\Psi^{\epsilon,\delta}))$ is compact in $N$.

Now let $\epsilon$ decrease to 0. 
Observe that $(\pi_2)_*(\log_\delta|z-z_0|\cdot R_{\theta_0})$ is defined by a p.s.h. function 
$\varphi_{\theta_0}^\delta$
on $N$ and
$$\Uc_{\epsilon,\delta}(\theta_0,z_0,w_0)=\int \varphi_{\theta_0}^\delta(w)\psi_\epsilon(w-w_0)\lambda(w-w_0).$$
Since $\psi$ is radial and $\varphi_{\theta_0}^\delta$ is p.s.h., the submean inequality implies that 
$\Uc_{\epsilon,\delta}$ decreases to a p.s.h. function $\Uc_\delta$
on $\Delta\times \C\times N$. The definition of $\Phi^{\epsilon,\delta}_{z_0,w_0}$ and slicing theory imply that
$$\Uc_\delta(\theta_0,z_0,w_0)= \big\langle\langle R_{\theta_0}, \pi_2, w_0\rangle,\log_\delta|z-z_0|\big\rangle.$$

Recall that $\langle R_\theta, \pi_2, w_0\rangle$ is a probability measure. When $\delta$ decreases to 0,
$\Uc_\delta$ decreases to the p.s.h. function 
$$\Uc(\theta_0,z_0,w_0) := \big\langle\langle R_{\theta_0}, \pi_2, w_0\rangle,\log|z-z_0|\big\rangle.$$
The last formula says that for every fixed $\theta_0$, $\Uc(\theta_0,\cdot,\cdot)$ defines a potential of
$R_{\theta_0}$. 
In particular,  the restriction of $\Uc$ to $\{\theta_0\}\times \C\times N$ is pluriharmonic
outside the support of $R_{\theta_0}$. 
Recall that the projection of  $\cup\supp(R_\theta)$ on $M$ is relatively compact in $M$.
On the other hand, for $|z_0|$ large enough, $\log_\delta|z-z_0|$ is pluriharmonic for $z\in M$.
Then, it is easy to check that $\Uc_{\epsilon,\delta}$ and $\Uc$ are pluriharmonic
for $|z_0|$ large enough. Now, by
Hartogs extension theorem, $\Uc$ is pluriharmonic 
near  $\Delta\times\partial_v D$ and then $\Rc:=\ddc\Uc$ vanishes near $\Delta\times \partial_v D$. 
It follows that the slices of $\Rc$, which
are equal to $R_\theta$, define a structural disc in $\Cc_v^1(D)$. 
\end{proof}

\section{Intersection of currents}
\label{intersection_currents}

In this section, we define the intersection (wedge product) $R\wedge S$ of a vertical positive closed current 
$R\in \Cc_v(D)$ 
and a horizontal
positive closed current $S\in\Cc_h(D)$. 
When one of these currents, for example $R$, 
has bidegree $(1,1)$, using a regularization, 
the reader can verify that our
definition coincides with the classical definition $R\wedge S:=\ddc(uS)$ where $u$ is a potential of $R$.
The current $uS$ is well defined since, by Oka's inequality 
\cite[Prop. 3.1]{FornaessSibony2}, 
$u$ is integrable with respect to the trace measure of $S$. This case is very simple since the mass of $uS$
on a compact set can be estimated using Stokes' theorem and the geometry of the supports of  $R$ and $S$.

\begin{theorem} \label{th_intersection_currents}
Let $R$ be a current in $\Cc_v(D)$ and $S$ be a current in $\Cc_h(D)$. Then 
$R\wedge S$ is defined such that for every p.s.h. function $\varphi$ on $D$
$$\langle R\wedge S,\varphi\rangle = \limsup_{R'\rightarrow R\atop S'\rightarrow S} \langle R'\wedge S',\varphi\rangle$$
where $R'\in \Cc_v(D)$ and $S'\in\Cc_v(D)$ are smooth with supports converging in the Hausdorff sense to 
those of $R$ and $S$. 
The value of $\langle R\wedge S,\varphi\rangle$ depends linearly on $R$, $S$ and $\varphi$.
The wedge product $R\wedge S$ is a positive measure of mass  $\|R\|_v\|S\|_h$ and supported 
in $\supp(R)\cap\supp(S)$. 
\end{theorem}

In the previous theorem, the convexity of $D$ is not necessary and
we can take $R'$ and $S'$
such that $\supp(R')\cap\supp(S')$ is contained in a fixed
neighbourhood $W\Subset D$ of $\supp(R)\cap\supp(S)$ (see Propositions
\ref{intersection_test} and \ref{intersection_continuity}).

Choose $M'$ and $N''$ such that $R\in\Cc_v(D')$ and $S\in\Cc_h(D'')$.
We can assume that $R$ and $S$ are normalized. We will construct explicitly the probability measure
$R\wedge S$. We first prove the following lemma.

\begin{lemma} \label{mass_of_intersection} 
Assume that 
$R$ is a continuous form. Then, $R\wedge S$ is a probability measure.
\end{lemma}
\begin{proof}
By regularization of currents, we can assume that $S$ is smooth. Let $\Rc$ be the structural 
disc associated to $R$ which was constructed
in Section \ref{geometry_currents}. 
The current 
$\Rc':=\Rc\wedge \Pi_D^*(S)$ is positive closed and of bidimension $(1,1)$. Moreover, the restriction of
$\Pi_V$ to $\supp(\Rc')$ is proper. Hence, $(\Pi_V)_*(\Rc')$, which is positive closed  and of bidimension $(1,1)$, 
is defined by a constant function $c$ 
on $V$. It follows that $\|R_\theta\wedge S\|=c$ for almost every $\theta$. By Lemmas \ref{disc1} and \ref{disc2},
$R_\theta$ depends continuously on $\theta$. Then $\|R_\theta\wedge S\|=c$ for every $\theta$.
From the explicit form of $R_0$, we get  
$c=\|R_0\wedge S\|=\|S\|_h=1$. For $\theta=1$, we get
$\|R\wedge S\|=c=1$.
\end{proof}

\begin{proposition} \label{compactness} 
Let $K$ be a compact subset of $M$. Then, the set of currents $R\in\Cc^1_v(D)$
with support in $K\times N$, is compact for the weak topology on currents.
\end{proposition}
\begin{proof} Let $L$ be a compact subset of $D$. Let $S\in\Cc^1_h(D)$ be a normalized smooth form,
strictly positive on $L$ (see lemma \ref{lemma_density}). Lemma \ref{mass_of_intersection} 
implies that $\langle R,S\rangle =1$. Hence, the mass of $R$
on $L$ is bounded from above by a constant independent of $R$. The proposition follows.
\end{proof}

Consider a function $\varphi$ continuous and p.s.h. in a neighbourhood $W'$ of $\supp(R)\cap \supp(S)$ in $D$.
Let $W$ be another neighbourhood of $\supp(R)\cap \supp(S)$ such that $W\Subset W'$.  
Consider smooth forms 
$R_n\in\Cc^1_v(D')$ and $S_n\in \Cc^1_h(D'')$ such that $R_n\rightarrow R$, $S_n\rightarrow S$,
$\supp(R_n)\cap\supp(S_n)\subset W$ and 
$\langle R_n\wedge S_n,\varphi\rangle$ converge to a constant $m_\varphi$. Assume that 
$m_\varphi$ is the maximal constant that we can obtain
in this way. It follows from Lemma \ref{mass_of_intersection} that $m_\varphi$ is finite.

Let $R_\theta$, $\theta\in V$,  
be the currents of the structural disc in $\Cc^1_v(D'')$ 
associated to $R$ that we constructed in Section \ref{geometry_currents}.
Recall that $R_1=R$.
We construct in the same way the horizontal currents $S_{\theta'}$, $\theta'\in V$, with $S_1=S$. They
define a structural disc in $\Cc^1_h(D')$. Observe that 
when $\theta,\theta'\rightarrow 1$, we have $\supp(R_\theta)\rightarrow\supp(R)$ and
$\supp(S_{\theta'})\rightarrow \supp(S)$. In particular, $\supp(R_\theta)\cap\supp(S_{\theta'})\subset W$
when $\theta$ and $\theta'$ are close to $1$.

\begin{proposition} \label{intersection_test}
We have 
$$m_\varphi=\limsup_{\theta\rightarrow 1} \langle R_\theta\wedge S,\varphi\rangle
=\limsup_{\theta\rightarrow 1}\langle R\wedge S_\theta,\varphi\rangle =\limsup_{\theta,\theta'\rightarrow
1}\langle R_\theta\wedge S_{\theta'},\varphi\rangle.$$
Moreover, $m_\varphi$ does not depend on $W'$,  $W$, and
depends linearly on $\varphi$, $R$, $S$. 
\end{proposition}
\begin{proof} 
Define 
$\psi(\theta,\theta'):=\langle R_\theta\wedge
S_{\theta'},\varphi\rangle$. By Lemma \ref{disc1} and Remark \ref{rk_nonconvex}, there exists a 
small neighbourhood $U$ of $(1,1)$ in $V^2$ such that
$\psi$ is defined and continuous
on $U\setminus(1,1)$. 
Lemma \ref{mass_of_intersection} shows that $\psi$ is bounded. 
We first show that $\psi$ is p.s.h. on
$U':=\{(\theta,\theta')\in U,\ \theta\not=1,\ \theta'\not=1\}$. This
allows us to extend $\psi$
to a p.s.h. function on $U$ with
$$\psi(1,1):=\limsup_{\theta,\theta'\rightarrow1}\psi(\theta,\theta') =
\limsup_{\theta,\theta'\rightarrow1} 
\langle R_\theta\wedge S_{\theta'},\varphi\rangle.$$

Let $\Rc$ and $\Sc$ be currents as in Section \ref{geometry_currents} whose slices are $R_\theta$ 
and $S_{\theta'}$. These currents are smooth for $\theta\not=1$ and $\theta'\not=1$. It follows
that the form 
$$\widetilde\Rc(\theta,\theta',z,w):=\varphi(z,w)\Rc(\theta,z,w)\wedge \Sc(\theta',z,w)$$ 
is continuous on $U'\times D$.
We also have 
$\ddc\widetilde\Rc\geq 0$ and the projection of $\supp(\widetilde\Rc)$ on $U'$ is proper.
As in Theorem \ref{th_slices}, we obtain $\psi$ as the push-forward of $\widetilde \Rc$ on $U'$. Hence
$\psi$ is p.s.h. 

Define $m'_\varphi:=\psi(1,1)$. 
We first prove that $m'_\varphi=m_\varphi$. This implies that $m_\varphi$ depends linearly on $\varphi$,
$R$ and $S$ since $\psi$ depends linearly on $\varphi$, $R$ and $S$.
The current $R_\theta$ is a priori not defined on $D$ but 
it is a vertical current on a domain $D_\theta$ with $D_\theta\rightarrow D$ when $\theta\rightarrow1$.
The current $S_{\theta'}$ satisfies the same properties.
Hence, by definition of $m_\varphi$, we have $m'_\varphi\leq m_\varphi$.

We define the structural discs $(R_{n,\theta})$ and 
$(S_{n,\theta})$ associated to $R_n$ and $S_n$ as in Section \ref{geometry_currents} with $R_{n,1}=R_n$ and 
$S_{n,1}=S_n$. Recall that $R_{n,\theta}$ and $S_{n,\theta}$ are smooth currents 
when $\theta\not =1$.
By Lemma \ref{disc1}, the bounded sequence of continuous p.s.h. functions 
$\psi_n(\theta,\theta'):=\langle R_{n,\theta}\wedge S_{n,\theta'},\varphi\rangle$ converges
to $\psi$ on $U\setminus(1,1)$. It follows that $\psi_n\rightarrow \psi$ in $\Loneloc(U)$. 
By Hartogs lemma, 
$$m'_\varphi=\psi(1,1)\geq\limsup_{n\rightarrow\infty} \psi_n(1,1) =
\limsup_{n\rightarrow\infty} \langle R_n\wedge S_n,\varphi\rangle =m_\varphi.$$
Hence $m_\varphi'=m_\varphi$. 

Since p.s.h. functions on $U$ are decreasing limits of 
 smooth p.s.h. functions, their restrictions to
$V\times \{1\}$ are subharmonic functions. It follows that 
$$\limsup_{\theta\rightarrow1}\langle R_\theta\wedge S,\varphi\rangle = \limsup_{\theta\rightarrow1}
\psi(\theta,1)=\psi(1,1)=m_\varphi.$$
We prove in the same way that $\limsup \langle R\wedge S_\theta,\varphi\rangle =m_\varphi$. 
\end{proof}

\noindent
{\bf End of the proof of Theorem \ref{th_intersection_currents}.} 
For functions $\varphi$ continuous p.s.h. on a neighbourhood of $\supp(R)\cap\supp(S)$, define 
$$\langle R\wedge S,\varphi\rangle :=m_\varphi.$$
Since smooth functions on neighbourhoods of $\supp(R)\cap\supp(S)$ 
can be written as differences of continuous p.s.h. functions, we can extend the definition to smooth functions.

Proposition \ref{intersection_test} shows that the current $R\wedge S$ is supported 
in $\supp(R)\cap \supp(S)$.
It is clear that the definition does not depend 
on coordinate systems of $M$, $N$. 
If $\varphi\leq \varphi'$ we have $m_\varphi\leq m_{\varphi'}$. Then $R\wedge S$ is a positive measure. When
$\varphi=1$, Lemma \ref{mass_of_intersection} implies that 
$m_\varphi=1$. Hence $R\wedge S$ is a probability measure.
\hfill $\square$

\begin{proposition} \label{intersection_continuity}
Let $R$, $S$, $\varphi$, $W'$ and $W$ be as above. 
Let $R_n\in\Cc_v(D')$ and $S_n\in \Cc_h(D'')$
such that $R_n\rightarrow R$, $S_n\rightarrow S$ and $\supp(R_n)\cap \supp(S_n)\subset W$. 
Then
$$\limsup_{n\rightarrow\infty} \langle R_n\wedge S_n,\varphi\rangle
\leq \langle R\wedge S,\varphi\rangle.$$ 
The measures $R_n\wedge S_n$ converge to $R\wedge S$ 
if and only if $\langle R_n\wedge S_n,\varphi\rangle \rightarrow \langle R\wedge S,\varphi\rangle$
for one function $\varphi$ strictly p.s.h. on $W'$.  
In particular, there exists $(\theta_n)\subset V\setminus\{1\}$ converging to $1$
such that 
$$R_{\theta_n}\wedge S\rightarrow R\wedge S, \quad R\wedge S_{\theta_n}
\rightarrow R\wedge S \quad \mbox{and} \quad  R_{\theta_n}\wedge S_{\theta_n}\rightarrow R\wedge S.$$
More generally, if $(\theta,\theta')\rightarrow(1,1)$ in the plurifine topology, then 
$R_\theta\wedge S_{\theta'}\rightarrow R\wedge S$. 
\end{proposition}
\begin{proof}
The first inequality follows from the definition of $m_\varphi$.
Now assume that $\varphi$ is strictly p.s.h. on $W'$. Let $\phi$ be a real smooth function with support in $W'$. 
If $A>0$ is large enough then $\varphi^\pm:=A\varphi\pm\phi$ are p.s.h. on $W'$. Then
$\limsup\langle R_n\wedge S_n,\varphi^\pm\rangle\leq \langle R\wedge S,\varphi^\pm\rangle$.
When $\langle R_n\wedge S_n,\varphi\rangle \rightarrow \langle R\wedge S,\varphi\rangle$, we deduce easily
that $\langle R_n\wedge S_n,\phi\rangle \rightarrow
\langle R\wedge S,\phi\rangle$. It follows that $R_n\wedge S_n
\rightarrow R\wedge S$.  

The functions $\psi(\cdot,1)$, $\psi(1,\cdot)$ and $\psi(\cdot,\cdot)$ associated to $\varphi$ are subharmonic or 
p.s.h. 
Then there exists $(\theta_n)\rightarrow1$
such that $\psi(\theta_n,1)$, $\psi(1,\theta_n)$ and
$\psi(\theta_n,\theta_n)$ converge to $\psi(1,1)$. Hence 
$\langle R_{\theta_n}\wedge S,\varphi\rangle=\psi(\theta_n,1)$
converge to $\langle R\wedge S,\varphi\rangle=\psi(1,1)$. 
It follows that $R_{\theta_n}\wedge S\rightarrow R\wedge S$. 
Other convergences are obtained in the same way. 
If $(\theta,\theta')\rightarrow(1,1)$ in the plurifine topology (i.e. the coarsest topology which makes 
p.s.h. functions continuous), we get 
$R_\theta\wedge S_{\theta'}\rightarrow R\wedge S$. 
\end{proof}

\begin{remarks} \label{intersection_rk}
\rm
{\bf a.} Proposition \ref{intersection_continuity} and 
Lemma \ref{disc2} 
imply that when $R$ or $S$ is continuous,
our definition of $R\wedge S$ coincides with the usual one. 

{\bf b.} When $\varphi$ is a uniform limit of continuous functions 
p.s.h. on  neighbourhoods of $\supp(R)\cap \supp(S)$,
we can apply 
Proposition \ref{intersection_continuity} and get $\limsup \langle R_n\wedge S_n,\varphi\rangle
 \leq \langle R\wedge S,\varphi\rangle$.
Hence, if there exists a compact set $K\subset D$ containing $\supp(R)\cap\supp(S)$ 
such that continuous p.s.h. functions 
on neighbourhoods of $K$ are dense in $\Cc^0(K)$, 
then $R_n\wedge S_n\rightarrow R\wedge S$
provided that $R_n\rightarrow R$, $S_n\rightarrow S$ and 
$\supp(R_n)\cap\supp(S_n)\rightarrow K$.   In particular, this holds
when $K$ is totally disconnected.
In the last case, continuous functions on $K$ can be approximated by
functions locally constant in neighbourhoods of $K$. 
\end{remarks}

Let $\lambda_\epsilon$ denote the Lebesgue measure on the disc of center $1$ and of radius $\epsilon$ 
normalized
by $\|\lambda_\epsilon\|=1$. Since the function $\psi$ in Proposition \ref{intersection_test} is p.s.h. we have
\begin{eqnarray*}
\psi(1,1) & = & \lim_{\epsilon\rightarrow 0} \int \psi(\theta,1)\d\lambda_\epsilon(\theta)
=\lim_{\epsilon\rightarrow 0} \int \psi(1,\theta)\d\lambda_\epsilon(\theta)\\
& = & \lim_{\epsilon\rightarrow 0}\int \psi(\theta,\theta')\d\lambda_\epsilon(\theta)\d\lambda_\epsilon(\theta').
\end{eqnarray*}
We define the vertical and horizontal currents $R^{(\epsilon)}$ and $S^{(\epsilon)}$ by
$$R^{(\epsilon)}:=\int R_\theta\d\lambda_\epsilon(\theta)  \mbox{ \ \ and \ \ } 
S^{(\epsilon)}:=\int S_\theta\d\lambda_\epsilon(\theta)$$
and deduce from the previous relations that
$$\psi(1,1)=\lim_{\epsilon\rightarrow 0} \langle R^{(\epsilon)}\wedge S,\varphi\rangle=
\lim_{\epsilon\rightarrow 0} \langle R\wedge S^{(\epsilon)},\varphi\rangle
= \lim_{\epsilon\rightarrow 0} \langle R^{(\epsilon)}\wedge S^{(\epsilon)},\varphi\rangle.$$
This and Proposition \ref{intersection_continuity}
imply the following result which can be considered as a
``less abstract'' definition of $R\wedge S$.
\begin{proposition} \label{intersection_explicite_construction}
Let $R$, $S$, $R^{(\epsilon)}$ and $S^{(\epsilon)}$ be as above. Then
$$R\wedge S=\lim_{\epsilon\rightarrow 0} R^{(\epsilon)}\wedge S =
\lim_{\epsilon\rightarrow 0} R\wedge S^{(\epsilon)} =\lim_{\epsilon\rightarrow 0}
R^{(\epsilon)}\wedge S^{(\epsilon)}.$$
\end{proposition}

\begin{remarks} \label{intersection_construction_rk}
\rm
It follows from the definition of $R\wedge S$ and  
from Proposition \ref{intersection_explicite_construction},  that for $\varphi$ p.s.h. on $W'$
$$\langle R\wedge S,\varphi\rangle = \limsup \langle R'\wedge S',\varphi\rangle$$
where the limit is taken over all currents, not 
necessarily smooth,  $R'\rightarrow R$ and
$S'\rightarrow S$ with $\supp(R')\cap\supp(S')\subset W$. 

Let $(R'_\theta)$ (resp. $(S'_{\theta'})$) be an {\it arbitrary}  
structural variety in $\Cc^1_v(D)$ (rep. in $\Cc^1_h(D)$).
Let $\varphi$ be a bounded p.s.h. function on $D$. Then one can prove as in Theorem \ref{th_slices} 
and Proposition \ref{intersection_test}
that the function 
$\lambda(\theta,\theta'):=\langle R'_\theta\wedge S'_{\theta'},\varphi\rangle$ is p.s.h. and 
$(\theta,\theta')\mapsto R_\theta'\wedge S_{\theta'}'$ is continuous for the plurifine topology.
\end{remarks}

\section{Horizontal-like maps}
\label{horiz_map}

In general, a 
horizontal-like map $f$ on $D$ is not defined on the whole domain $D$ but only on 
a vertical subset $f^{-1}(D)$ of $D$. It takes values in a horizontal subset 
$f(D)$ of $D$.  
We define these maps using their graphs as follows 
(see \cite{Dujardin, DinhDujardinSibony}). Let $\pr_1$ and $\pr_2$ be the canonical 
projections of $D\times D$ on its factors. From now on we assume
always that $D$ is convex.

\begin{definition} \rm \label{def_horiz_map}
A {\it horizontal-like map} $f$ on $D$ is a holomorphic map with graph 
$\Gamma$ such that
\begin{enumerate}
\item $\Gamma$ is an irreducible submanifold of $D\times D$.
\item $\pr_{1|\Gamma}$ is injective; $\pr_{2|\Gamma}$ has finite fibers.
\item $\overline\Gamma$ does not intersect 
$\overline{\partial_v D} \times \overline D$ nor 
$\overline D \times \overline{\partial_h D}$.
\end{enumerate} 
\end{definition}
The map $f$ is defined on $f^{-1}(D):=\pr_1(\Gamma)$ 
and its image is equal to $f(D):=\pr_2(\Gamma)$
(if we assume only that $\pi_{1|\Gamma}$ has finite fibers, we obtain a horizontal-like correspondence).
Observe that there exist open sets 
$M'\Subset M$ and $N''\Subset N$ such that
$f^{-1}(D)\subset D':=M'\times N$ and $f(D)\subset D'':=M\times
N''$. We have $\Gamma\subset D'\times D''$. This property
characterizes horizontal-like maps.
Since $\Gamma$ is a submanifold of $D\times D$, 
when $x$ converges to $\partial f^{-1}(D)\cap D$, 
$f(x)$ converges to $\partial_vD$.
When $y$ converges to $ \partial f(D)\cap D$, $f^{-1}(y)$ converges to 
$\partial_hD$.
So, the vertical (resp. horizontal) part of $\partial f^{-1}(D)$ is sent to the
vertical (resp. horizontal) part of $\partial f(D)$.
If $g$ is another horizontal-like map on $D$, $f\circ g$ is also a horizontal-like map.
When $p=k$, we obtain the polynomial-like maps studied in 
\cite{DinhSibony1}.

If $\pr_{2|\Gamma}$ is injective, we say that $f$ is {\it invertible}. 
In this case, up to a coordinate change, 
$f^{-1}:\pr_2(\Gamma)\rightarrow\pr_1(\Gamma)$ is a horizontal-like map.
When $k=2$ and $p=1$,  we obtain the H{\'e}non-like maps which are studied in 
\cite{Dujardin, DinhDujardinSibony}. 
In order to simplify the paper, we consider only invertible horizontal-like maps. The results in Sections 
\ref{horiz_map}, \ref{random_iteration} and 
\ref{section_green_current} hold for non-invertible maps, but for the construction of $T_+$, we need to define 
inverse images of positive closed currents by open holomorphic maps, see also 
\cite{DinhSibony4}.

The operator $f_*=(\pr_{2|\Gamma})_*\circ (\pr_{1|\Gamma})^*$ 
acts continuously on horizontal currents. 
If $S$ is a horizontal current (form), so is 
$f_*(S)$.
The operator $f^*=(\pr_{1|\Gamma})_*\circ (\pr_{2|\Gamma})^*$ 
acts continuously on vertical currents. 
If $R$ is a vertical current (form), so is $f^*(R)$. 
We have the following proposition for positive closed currents.

\begin{proposition} \label{operators_on_currents} The operator $f_*:\Cc_h(D')\rightarrow \Cc_h(D'')$
is well defined and continuous. Moreover, there exists an integer $d\geq 1$ such that
$\|f_*(S)\|_h=d\|S\|_h$ for every $S\in \Cc_h(D')$. 
The operator $f^*:\Cc_v(D'')\rightarrow \Cc_v(D')$
is well defined and continuous. If $R$ belongs to $\Cc_v(D'')$, we have 
$\|f^*(R)\|_v=d\|R\|_v$.
\end{proposition}

\begin{proof} Using Definition \ref{def_horiz_map}, one can check that $f^*$ and $f_*$ are well
defined.

Let $R$ be a current in $\Cc^1_v(D'')$. We want to compute the slice mass
of $f^*(R)$. We can assume that $R$ is smooth. Let $S=[w=b]$ be the current 
of integration on the subspace $\{w=b\}$ with $b\in N$. Since $S$ is normalized, 
we have 
$$\|f^*(R)\|_v=\langle f^*(R), S \rangle =\langle R, f_*(S)\rangle.$$
The current $f_*(S)$ is defined by a horizontal analytic subset of $D''$. Hence, it is a 
ramified covering of degree $d$ over $M$. We have $\|f_*(S)\|_h=d$. 
Theorem \ref{th_intersection_currents} implies that $\langle R, f_*(S)\rangle =d$. 
Hence, $\|f^*(R)\|_v=d$. 

If $S$ is an arbitrary current in $\Cc^1_h(D')$, then 
Theorem \ref{th_intersection_currents}
implies that
$$\|f_*(S)\|_h=\langle f_*(S),R\rangle =\langle S, f^*(R)\rangle =d.$$ 
\end{proof}

The integer $d$ in Proposition \ref{operators_on_currents} is called
the {\it (main) dynamical degree} of $f$. Define 
$\L_v:=\frac{1}{d} f^*$ and 
$\L_h:=\frac{1}{d} f_*$. Using Ces{\`a}ro means, one can easily construct  
a current $T_+\in \Cc_v^1(D)$ such that $\L_v(T_+)=T_+$. 
A priori such $T_+$ is not unique.
Our aim is to construct such a current $T_+$ with a 
good convergence theorem and some extremality properties. This allows us to construct 
an interesting invariant 
measure.
The following diagram is one of the main objects 
we consider:
\begin{eqnarray}
\xymatrix{
          V \ar[r]^-{\tau}  & 
\Cc^1_v(D^*) \ar@(dl,dr)_{\L_v}[]   \ar[r]^-{\Lambda_\Phi} & \R
} \label{diagram2}
\end{eqnarray}

\begin{example}\rm \label{ex_regular_maps}
Let $f$ be a polynomial automorphism of $\C^k$. Denote also by $f$ its meromorphic extension to $\mathbb{P}^k$.
Let $(z_1,\ldots, z_k)$ be the coordinates of $\C^k$ and $[z_1:\cdots:z_k]$ be homogeneous coordinates of
the hyperplane at infinity $L$. 
Assume that the indeterminacy set $I_+$ of $f$ is the subspace $\{z_1=\cdots =z_p =0\}$
of $L$ and the indeterminacy set $I_-$ of $f^{-1}$ is the subspace $\{z_{p+1}=\cdots=z_k=0\}$ of $L$.
This map is regular in the sense of \cite{Sibony}; that is $I_+\cap I_-=\varnothing$ (see also \cite{DinhSibony5}). 

If $M$ and $N$ are the balls of center $0$ and of radius $r$ in $\C^p$ and $\C^{k-p}$, 
then $f^{n_0}$ defines a horizontal-like map in $D=M\times N$ when $r$ and $n_0$ 
are big enough. 
This follows from the description of Julia sets
of $f$ and $f^{-1}$ in \cite{Sibony}. 

Observe that every small pertubation of $f^{n_0}$ on $D$ is still horizontal-like. 
One can construct such a 
map which admits both attractive and repelling fixed points \cite{Dujardin}. The map is not conjugated
to a polynomial automorphism since polynomial automorphisms have constant jacobian and hence cannot 
have such fixed points. 
\end{example}

\begin{example}\rm \label{ex_product_maps}
Let $f_i$ be horizontal-like maps on $D_i=M_i\times N_i$. Define $D=D_1\times D_2$ and 
the product map $f(x_1,x_2):=(f_1(x_1),f_2(x_2))$. Up to a coordinate change, we can identify $D$
to $M\times N$, with $M=M_1\times M_2$ and $N=N_1\times N_2$. Then, one can check easily that $f$ defines a
horizontal-like map on $D$.

When $M_1=N_2$ and $N_1=M_2$, let $\Delta$ denote the diagonal of $D$. Then, $\Delta$ is not a horizontal set but 
$f(\Delta)$ is horizontal.

We will see in Section \ref{section_equilibrium_measure} 
that this simple example can be used to linearize some problems.
\end{example}

\section{Random iteration}
\label{random_iteration}

Let $(f_n)$ be a sequence of invertible
horizontal-like maps on $D$ of dynamical degrees $d_n$.
Define $\L_{v,n}:=\frac{1}{d_n}f_n^*$ and $\L_{h,n}:=\frac{1}{d_n}(f_n)_*$.
Assume there exist open sets $M'\Subset  M$ and $N''\Subset N$ such that 
$f_n^{-1}(D)\subset D':=M'\times N$  and $f_n(D)\subset D'':=M\times N''$ 
for every $n$. 
Define the {\it filled Julia set} associated to $(f_n)$ as 
$$\Kc_+:=\bigcap_{n\geq 1} f_1^{-1}\circ\cdots \circ f_n^{-1}(D)
=\bigcap_{n\geq 1} f_1^{-1}\circ\cdots \circ f_n^{-1}(\overline D').$$ 
This is a vertical closed subset of $D'$.

\begin{theorem} \label{random_green_current} Let $(R_n)\subset \Cc^1_v(D')$ be a uniformly bounded family of 
forms.
Then, the sequence $\L_{v,1}\ldots\L_{v,n}(R_n)$ converges 
weakly to a current $T_+\in\Cc^1_v(D')$
supported in $\partial \Kc_+$. Moreover, $T_+$ is independent of $(R_n)$.
\end{theorem}   

We say that $T_+$ is the {\it Green current} associated to the sequence $(f_n)$. 
We say that $(R_n)$ is {\it uniformly bounded} if the coefficients of $R_n$ are uniformly bounded.
Observe that $\L_{v,n}$ is ``distance decreasing'' for the Kobayashi pseudo-distance on $\Cc_v^1(\overline D'')$. 
However, the fact that it is not a distance makes the convergence questions more delicate. 
We first prove the following proposition.

\begin{proposition} \label{random_test}
Let $\Phi$ be a real continuous horizontal $(k-p,k-p)$-form with $\ddc\Phi\geq 0$.
There exists a constant $M_\Phi$ such that
if $R_n$ are currents in $\Cc^1_v(D)$, then 
$\limsup\langle \L_{v,1}\ldots \L_{v,n} (R_n),\Phi\rangle\leq M_\Phi$; 
if $R_n$ are as in Theorem \ref{random_green_current}, then
$\lim \langle \L_{v,1}\ldots \L_{v,n} (R_n),\Phi\rangle=M_\Phi$.
\end{proposition} 

\begin{proof}
By regularization, we can assume that $R_n$ are smooth.
Observe that by Theorem \ref{th_intersection_currents} if $\Phi$ is positive and closed then 
$\langle \L_{v,1}\ldots \L_{v,n} (R_n),\Phi\rangle=\|\Phi\|_h$.
So in this case the convergence is clear.
If we add to $\Phi$ a form in $\Cc_h(D)$, we can assume that
$\Phi$ is positive on $D'$. We can also assume that $\Phi$ is smaller than a smooth 
form in $\Cc^1_h(D')$. It follows from Proposition \ref{operators_on_currents} that each form 
$\L_{h,n}\ldots\L_{h,1} (\Phi)$
is positive and bounded from above by a current in $\Cc^1_h(D'')$,
depending on $n$.

Let $(\widetilde R_{i_n}')$ be a sequence of continuous 
forms in $\Cc^1_v(D'')$ with $(i_n)$ 
a sequence of integers, $i_n>n$, such that
$\langle \L_{v,1}\ldots\L_{v,i_n} (\widetilde R_{i_n}'),\Phi\rangle$ converge
to a real number $M_\Phi$. 
We choose $(i_n)$ and $(\widetilde R_{i_n}')$ so that $M_\Phi$ is the maximal value
that we can obtain in this way.
Since  $\L_{h,n}\ldots\L_{h,1} (\Phi)$ are bounded by normalized currents, Theorem \ref{th_intersection_currents}
implies that $M_\Phi$ is finite.
Hence, $M_\Phi$ satisfies the inequality in Proposition \ref{random_test}.

Define $\widetilde R_n:=\L_{v,n+1}\cdots\L_{v,i_n} (\widetilde R_{i_n}')$. We have 
$\widetilde R_n\in \Cc^1_v(D')$ and 
$\langle \L_{v,1}\ldots\L_{v,n} (\widetilde R_n),\Phi\rangle\rightarrow M_\Phi$.
We will use the 
structural discs $(\widetilde R_{n,\theta})$ of $\Cc^1_v(D'')$ 
constructed in Section \ref{geometry_currents} (see also (\ref{diagram1}) and (\ref{diagram2})) associated 
to $\widetilde R_n$ in order to prove that the convergence holds when
$R_n$ is replaced by $\widetilde R_{n, 0}$.

Theorem \ref{th_slices} allows us to define continuous subharmonic functions on $V$ by
$$\varphi_n(\theta):=\langle \L_{v,1}\ldots\L_{v,n} (\widetilde R_{n,\theta}),\Phi\rangle
=\langle \widetilde R_{n,\theta}, \L_{h,n}\ldots \L_{h,1}(\Phi)\rangle.$$
Since $\varphi_n(1)$ tends to the maximal value 
$M_\Phi$, Hartogs lemma \cite{Hormander} and the maximum principle imply that
$\varphi_n\rightarrow M_\Phi$ in $\Loneloc(V)$. 

On the other hand, since each 
$\L_{h,n}\ldots \L_{h,1}(\Phi)$ is bounded by a 
current in $\Cc^1_h(D)$, Lemma \ref{disc1}
implies that $|\varphi_n(\theta)-\varphi_n(0)|\leq c|\theta|$ 
for $|\theta|\leq r$.
Hence, $\varphi_n(0)$ converge to $M_\Phi$.
Since $\widetilde R_0:=\widetilde R_{n,0}$ is independent of $n$, 
we obtain that $\langle \L_{v,1}\ldots\L_{v,n} (\widetilde R_0),
\Phi\rangle\rightarrow M_\Phi$.

Now assume that $R_n$ satisfy the hypothesis of Theorem \ref{random_green_current}.
If we replace $M'$ by 
a bigger domain, we can assume that there
exists an open set $M''\Subset M'$ such that $f_n^{-1}(D)\subset M''\times N$ and $\supp(R_n)\subset M''\times N$.
Then, we can find a continuous form $R\in\Cc^1_v(D')$ and $c>0$ 
such that $R_n\leq cR$ for every $n$.

Define the currents $R_\theta$ associated to $R$ as in Section \ref{geometry_currents} and  
$$\psi_n(\theta):=\langle \L_{v,1}\ldots\L_{v,n} (R_\theta),\Phi\rangle
=\langle R_\theta, \L_{h,n}\ldots \L_{h,1}(\Phi)\rangle.$$
Recall that $R_0=\widetilde R_0$.
Since $\psi_n(0)=\varphi_n(0)\rightarrow M_\Phi$
and $\limsup \psi_n\leq M_\Phi$, we have 
$\psi_n\rightarrow M_\Phi$ in $\Loneloc(V)$. 
On the other hand, since $\L_{h,n}\ldots \L_{h,1}(\Phi)$ are bounded by 
currents in $\Cc^1_h(D'')$, Lemma \ref{disc2} implies that
$$\lim_{\theta\rightarrow1} 
\left(\sup_{n\geq 1} |\psi_n(\theta)-\psi_n(1)|\right) =0.$$
It follows that $\psi_n(1)\rightarrow M_\Phi$. We obtain that
$\langle \L_{v,1}\ldots\L_{v,n} (R),\Phi\rangle
\rightarrow M_\Phi$.

We turn to the general case.
Since $R_n$ and $cR -R_n$ belong to $\Cc_v(D')$, 
by definition of $M_\Phi$, we have 
\begin{eqnarray}
\limsup\langle \L_{v,1}\ldots \L_{v,n} (R_n),\Phi\rangle\leq M_\Phi \label{ineq_M_Phi}
\end{eqnarray}
and since $cR-R_n$ have slice mass $c-1$
\begin{eqnarray}
\limsup\langle \L_{v,1}\ldots \L_{v,n} (cR-R_n),\Phi\rangle\leq (c-1)M_\Phi. \label{ineq_M_Phi_bis}
\end{eqnarray}
We consider the sum of (\ref{ineq_M_Phi}) and (\ref{ineq_M_Phi_bis}) and deduce that these 
inequalities are in fact equalities. It follows that
$\lim\langle \L_{v,1}\ldots \L_{v,n} (R_n),\Phi\rangle= M_\Phi$.
\end{proof}

\begin{remark} \label{rk_random_test}
\rm 
Proposition \ref{random_test} still holds when $R_n$ are continuous forms and 
$\Phi$ is a non-smooth horizontal current such that $\ddc \Phi\geq 0$
and $-\Psi\leq \Phi\leq \Psi$ for some current $\Psi\in \Cc_h(D)$.
If $R_n$ are continuous then $\L_{v,1}\ldots \L_{v,n}(R_n)$ are continuous
and they act on currents of order 0, like $\Phi$. All the arguments in the above 
proof make sense. In this case, the convergence in Theorem \ref{random_green_current} is stronger than
the usual weak convergence.    
\end{remark}

\noindent
{\bf Proof of Theorem \ref{random_green_current}.} 
Since the maps $\Lambda_\Phi$ separate the currents 
in $\Cc_v(D)$, Proposition \ref{random_test} implies that 
$\L_{v,1}\ldots \L_{v,n} (R_n)$
converge to a current $T_+$ in $\Cc^1_v(D)$ which is defined
by  $\langle T_+, \Phi\rangle := M_\Phi$. This current is independent of $(R_n)$.

Now, we prove that $T_+$ is supported in $\partial \Kc_+$. It is clear that $\supp(T_+)\subset \Kc_+$.    
If $U\Subset \Kc_+$ is an open set, then 
$f_n\circ\cdots \circ f_1(U)\subset f_{n+1}^{-1}(D)\subset M''\times N$ for some $M''\Subset M$ and for every $n$. 
It follows that if $\supp(R_n)\subset (M'\setminus\overline M'')\times N$ we get
$\supp(T_+)\cap U=\varnothing$, since $T_+$ is independent of $R_n$.
\hfill $\square$
\\

The following corollary is a direct consequence of Proposition \ref{random_test}. It
gives an extremality property of $T_+$.

\begin{corollary} \label{extremality_random_current}
Let $(R_n)\subset \Cc^1_v(D)$.
Let $\Phi$ be a real continuous horizontal $(k-p,k-p)$-form such that $\ddc\Phi\geq 0$.
Then, every limit value $R$
of the sequence of currents $\L_{v,1}\ldots \L_{v,n} (R_n)$ satisfies 
$$\langle R,\Phi\rangle \leq \langle T_+,\Phi \rangle.$$ 
If $\ddc\Phi=0$, then $\langle R,\Phi\rangle =\langle T_+,\Phi\rangle$. 
\end{corollary}

\begin{corollary} \label{random_dominating_convergence} Let $(n_i)$ be an increasing sequence of integers and 
$R_{n_i}$, $R_{n_i}'$ be currents in $\Cc^1_v(D)$.
Assume that $\L_{v,1}\ldots \L_{v,n_i}(R_{n_i})$ converge to $T_+$ and that $R_{n_i}'\leq c R_{n_i}$
with $c>0$ independent of $n_i$. Then, $\L_{v,1}\ldots \L_{v,n_i}(R_{n_i}')$ converge also to $T_+$.
\end{corollary}
\begin{proof}
Let $\Phi$ be as above. Proposition \ref{random_test} implies that 
\begin{eqnarray}
\limsup \langle \L_{v,1}\ldots \L_{v,n_i}(R'_{n_i}),\Phi\rangle 
& \leq & \langle T_+,\Phi\rangle. \label{ineq_T_plus}
\end{eqnarray}
On the other hand, the currents $cR_{n_i}-R_{n_i}'$ belong to 
$\Cc_v(D)$ and have slice mass $c-1$.
Hence 
\begin{eqnarray}
\limsup \langle \L_{v,1}\ldots\L_{v,n_i}(cR_{n_i}-R'_{n_i}),\Phi\rangle 
\leq (c-1)\langle T_+,\Phi\rangle. \label{ineq_T_plus_bis}
\end{eqnarray}
By hypothesis,
\begin{eqnarray*}
\lim \langle \L_{v,1}\ldots \L_{v,n_i}(R_{n_i}),\Phi\rangle 
& = & \langle T_+,\Phi\rangle.
\end{eqnarray*}
We consider the sum of (\ref{ineq_T_plus}) and (\ref{ineq_T_plus_bis}) and 
deduce that 
$$\lim \langle \L_{v,1}\ldots\L_{v,n_i}(R_{n_i}'),\Phi\rangle =\langle T_+,\Phi\rangle.$$
The corollary follows.
\end{proof}

The following 
proposition allows us to check that $\lim\L_{v,1}\ldots \L_{v,n}(R_n)= T_+$ with only one test form.

\begin{proposition} \label{random_one_test_form}
Let $R_n$, $R$ and $\Phi$ be as in Corollary \ref{extremality_random_current}.
Assume that $\ddc\Phi$ is strictly positive on an open set $V$. 
If $\langle R,\Phi\rangle =\langle T_+,\Phi\rangle$, then $R=T_+$ on $V$. 
\end{proposition}
\begin{proof} Let $\Psi$ be a real test form with compact support in $V$. Let $A>0$ be a constant such that
$\ddc(A\Phi\pm\Psi)\geq 0$. Corollary \ref{extremality_random_current} implies that 
$\langle R,A\Phi\pm\Psi\rangle \leq \langle T_+, A\Phi\pm \Psi\rangle$. 
We deduce that $\langle R,\Psi\rangle = \langle T_+,\Psi\rangle$
if $\langle R,\Phi\rangle = \langle T_+,\Phi\rangle$.
Therefore, $R=T_+$ on $V$. 
\end{proof}

\begin{corollary} \label{random_varieties_convergence}
Let $(n_i)$ be an increasing sequence of integers. Then, there exist a subsequence
$(m_i)$ and a pluripolar set $\Ec_+\subset M$ such that, for every $a\in M\setminus \Ec_+$, we have
$$\L_{v,1}\cdots \L_{v,m_i}[z=a]\rightarrow T_+$$
where $[z=a]$ is the current of integration on the vertical analytic set $\{a\}\times N$.
\end{corollary}
\begin{proof}
Let $\Phi$ and $V$ be as above. Consider the 
locally uniformly bounded p.s.h. functions 
$\varphi_{n_i}(a):=\langle \L_{v,1}\ldots \L_{v,n_i}[z=a],\Phi\rangle$ 
(see Section \ref{geometry_currents}). 
By extracting a subsequence, we can assume that $\varphi_{n_i}$
converge in $\Loneloc(M)$ to a p.s.h. function $\varphi$. Proposition \ref{random_test} implies 
$\varphi\leq M_\Phi$.

Let $\nu$ be a smooth probability measure with compact support in $M'$. Consider 
the current $R:=\pi_1^*(\nu)$ in $\Cc^1_v(D')$. Since $R$ is smooth, Proposition \ref{random_test} 
implies
$$\int\varphi_{n_i}\d\nu = \langle \L_{v,1} \cdots \L_{v,n_i}(R),
\Phi\rangle \rightarrow M_\Phi.$$
It follows that $\varphi=M_\Phi$. Hence, there exists a subsequence $(m_i)\subset (n_i)$ and a pluripolar set 
$\Ec_+(\Phi)\subset M$ such that
$\varphi_{m_i}\rightarrow M_\Phi$ pointwise on $M\setminus\Ec_+(\Phi)$ \cite[Proposition 3.9.4]{DinhSibony1}. 
Proposition \ref{random_one_test_form} implies that
$\L_{v,1}\ldots \L_{v,m_i}[z=a]\rightarrow T_+$ on $V$ for $a\not\in \Ec_+(\Phi)$.

Consider a sequence of $(\Phi_n, V_n)$ such that $\cup_n V_n=D$. 
Extracting subsequences of $(m_i)$ gives
$\L_{v,1}\ldots \L_{v,m_i}[z=a]\rightarrow T_+$ on $D$ for 
$a\not\in \Ec_+:=\cup_n \Ec_+(\Phi_n)$.
\end{proof}

\begin{remark} \label{rk_random_varieties}
\rm 
Corollary \ref{random_varieties_convergence} implies that $T_+$ can be approximated by currents of 
integration on vertical manifolds with control of support.
When $p=1$, this holds for every current in $\Cc_v(D)$ 
\cite{DuvalSibony}. The problem is still
open for general currents of higher bidegree. 
\end{remark}

\section{Green currents}
\label{section_green_current}

In the rest of the paper, we study the dynamics of an invertible horizontal-like map. 
The following result is a direct consequence of Theorem \ref{random_green_current}. 

\begin{corollary} \label{green_current}
Let $f$ be an invertible 
horizontal-like map on $D$ of dynamical degree $d\geq 1$. 
Let $\Kc_+:=\cap_{n\geq 1} f^{-n}(D)$ be the filled Julia set of $f$. 
Let $(R_n)\subset \Cc^1_v(D')$ be a uniformly bounded family of 
forms. Then, 
$d^{-n}f^{n*}(R_n)$ converge weakly 
to a current $T_+\in\Cc^1_v(D')$ supported in $\partial \Kc_+$. 
Moreover, $T_+$ does not depend on $(R_n)$
and satisfies $f^*(T_+)=dT_+$.
\end{corollary}

We call $T_+$ the {\it Green current} of $f$. 
Corollary \ref{random_varieties_convergence} 
shows that $T_+$ is a limit value of $(d^{-n}f^{n*}[z=a])$ for $a\in M$ generic.
We construct in the same way
the Green current $T_-\in \Cc^1_h(D'')$ for $f^{-1}$. This current is supported in 
the boundary of $\Kc_-:=\cap_{n\geq 1} f^n(D)$ and satisfies
$f_*(T_-)=dT_-$. 
Now, we give some properties of the Green currents. 

Let $(R_n)$ be an arbitrary sequence of currents in $\Cc^1_v(D)$ and $\Phi$ be a
smooth real horizontal test form such that $\ddc\Phi\geq 0$. Corollary \ref{extremality_random_current} 
implies that every limit value
$R$ of $(d^{-n}f^{n*}(R_n))$ satisfies $\langle R,\Phi\rangle\leq \langle T_+,\Phi\rangle$.  
Proposition \ref{random_one_test_form} 
implies that if $\langle R,\Phi\rangle =\langle T_+,\Phi\rangle$, then $R=T_+$ in the open
set where $\ddc\Phi$ is strictly positive. We deduce from this the following corollary.

\begin{corollary} \label{one_test_form}
Let $T$ be a current in $\Cc^1_v(D)$ and $\Phi$ be a real horizontal continuous
form. Assume that $\ddc\Phi\geq 0$ on $D$ and 
$\ddc \Phi>0$ on a neighbourhood $W$ of $\Kc_+\cap \Kc_-$. Then,  
$d^{-n}f^{n*}(T)\rightarrow T_+$ if and only if 
$\langle d^{-n}f^{n*}(T),\Phi\rangle\rightarrow \langle T_+,\Phi\rangle$. 
\end{corollary}
\begin{proof} 
Assume that $\langle d^{-n}f^{n*}(T),\Phi\rangle\rightarrow \langle T_+,\Phi\rangle$. 
Hence, every limit value $R$ of $(d^{-n}f^{n*}(T))$ is equal to $T_+$ on $W$. 
For every $m\geq 0$, there exists a limit value $R'$ of $(d^{-n+m}(f^{n-m})^*(T))$ such that $R=d^{-m}f^{m*}(R')$.
We also have $R'=T_+$ on $W$. This implies $R=T_+$ on $f^{-m}(W)$. It follows that $R=T_+$ on
$\cup_{m\geq 0} f^{-m}(W)$ which is a neighbourhood 
of $\Kc_+$. Since both the currents $R$ and $T_+$ are supported
in $\Kc_+$, we have $R=T_+$.
\end{proof}

The following result is a
direct consequence of Corollary \ref{random_dominating_convergence}.

\begin{corollary} \label{dominating_convergence} 
Let $T$ be a current in $\Cc^1_v(D)$.
Assume there exist $c>0$, 
an increasing sequence $(n_i)$ and currents $T_{n_i}\in \Cc^1_v(D)$
such that $T_{n_i}\leq cT_+$ and 
$T=d^{-n_i} (f^{n_i})^*(T_{n_i})$. Then $T=T_+$.
In particular, $T_+$ is extremal in the cone of currents 
$T\in\Cc^1_v(D)$
satisfying $f^*(T)=dT$.
\end{corollary}

\begin{theorem} \label{non_closed_current}
Let $R$ be a real continuous 
vertical form of bidegree $(p,p)$ not necessarily closed.
Then, $d^{-n}(f^n)^*(R)$ converge to $cT_+$ where $c:=\langle R,T_-\rangle$.
\end{theorem}
\begin{proof} 
We can write $R$ as a difference of positive forms (scale $D$ if necessary). 
Hence, we can assume that $R$ is positive and that $R\leq R'$ for  
a suitable
continuous form $R'\in \Cc_v(D)$. 
We can extract from $d^{-n}(f^n)^*(R)$ convergent subsequences. 
Corollary \ref{green_current} implies that every limit value is bounded by $\|R'\|_v T_+$. 

Let $(n_i)$
and $T$ such that $\lim d^{-n_i}(f^{n_i})^*(R)=T$. We have $T\leq \|R'\|_v T_+$.
Moreover, for every $m\geq 0$, we have $T=d^{-m}(f^m)^*(T')$ where $T'$ is a limit value of
$(d^{-n_i+m}(f^{n_i-m})^*(R))$.

Let $\Phi\in\Cc^1_h(D)$ be a continuous form. We have 
$$\langle T,\Phi\rangle = \lim \langle d^{-n_i}(f^{n_i})^*(R),\Phi\rangle
=\lim \langle R,d^{-n_i}(f^{n_i})_*(\Phi)\rangle =\langle R, T_-\rangle = c.$$
It follows that if $T$ were closed, it has slice mass $c$ (this also holds for $T'$). 
Hence, Corollary \ref{dominating_convergence} implies that it is sufficient to prove that 
$T$ is closed. 
We first prove that it is $\ddc$-closed.

\begin{lemma} \label{ddc_closed_limit} Let $T$ be a real vertical current of bidegree $(p,p)$ 
and of finite mass. 
Consider smooth forms
$\Phi\in \Cc^1_h(D)$.
Assume that 
$\langle T,\Phi\rangle$ does not depend on $\Phi$. Then $T$ is $\ddc$-closed.
\end{lemma}
\begin{proof} Consider a real smooth $(k-p-1,k-p-1)$-form $\alpha$ with compact support in $D$. 
Let $\Phi$ be a smooth form in $\Cc_h(D)$ strictly positive in 
a neighbourhood of $\supp(\alpha)$.
Write $\ddc\alpha=(A\Phi+\ddc\alpha) -A\Phi$. When $A$ is big enough, 
both $A\Phi+\ddc\alpha$ and $A\Phi$ are positive  
closed and have the same slice mass. By hypothesis, $\langle T,A\Phi+\ddc\alpha\rangle = \langle T,A\Phi\rangle$.
Hence, $\langle T,\ddc\alpha\rangle =0$ and $T$ is $\ddc$-closed. 
\end{proof}

Consider the product map $F(x_1,x_2)=(f(x_1), f(x_2))$ on $D^2$ as in Example \ref{ex_product_maps}. 
The same arguments
applied to $F$ and to $R\otimes R$ imply that $T\otimes T$ is $\ddc$-closed. It follows that $T$ is closed.
For this end, it suffices to compute $\ddc(T\otimes T)$.  
\end{proof}

\section{Equilibrium measure}
\label{section_equilibrium_measure}
The main result of this section is the following theorem.
\begin{theorem} \label{equilibrium_measure}
Let $f$ be an invertible horizontal-like map of dynamical degree $d$
on $D$. 
Let $(R_n)\subset \Cc^1_v(D')$ and $(S_n)\subset \Cc^1_h(D'')$ be uniformly 
bounded sequences of continuous 
forms. Then, $d^{-2n}(f^n)^*(R_n)\wedge (f^n)_*(S_n)$ converge
weakly to an invariant 
probability measure $\mu$ which does not depend on $(R_n)$ and $(S_n)$. Moreover,
$\mu$ is mixing and is supported on the boundary of the compact set $\Kc:=\cap_{n\in \Z} f^n(D)$.
\end{theorem}

We say that $\mu$ is the {\it equilibrium measure} of $f$. We will see that the  
convergence part of Theorem \ref{equilibrium_measure} is a consequence 
of Proposition \ref{random_test} and Remark \ref{rk_random_test} (see also Proposition 
\ref{stronger_convergence_measure} and Corollary \ref{convergence_on_current}).

Let $M_i$ and $N_i$ be copies of $M$ and $N$. Consider the domain
$$D^2=D\times D =(M_1\times N_1)\times (M_2\times N_2)\subset\C^{2k}$$ 
and the product map (see Example \ref{ex_product_maps})
$$F(z_1,w_1, z_2,w_2):=\big(f(z_1,w_1), f^{-1}(z_2,w_2)\big).$$ 
Using the coordinate change $(z_1,w_1,z_2,w_2)\mapsto (z_1,w_2,z_2,w_1)$, 
write
$$F(z_1,w_2,z_2,w_1)=\left(f_M(z_1,w_1), f_N^{-1}(z_2,w_2), f_M^{-1}(z_2,w_2), 
f_N(z_1,w_1)\right)$$ 
where $f=(f_M,f_N)$ and $f^{-1}=(f^{-1}_M,f^{-1}_N)$. 

Recall that, the coordinate change $(z,w)\mapsto (w,z)$ makes $f^{-1}$ a horizontal-like map.
One can check that
$F$ is an invertible  horizontal-like map of dynamical degree $d^2$ on $D^2\simeq 
(M_1\times N_2)\times (M_2\times N_1)$. The diagonal 
$$\Delta:=\{z_1=z_2, w_1=w_2\}$$
is not a horizontal set but $F(\Delta)$ is horizontal. If $\widetilde\varphi$ 
is a (positive) p.s.h. function on $\Delta$, then $\widetilde\varphi[\Delta]$ 
is a (positive) current such that $\ddc(\widetilde\varphi[\Delta])\geq 0$. 
Hence, we can
apply Proposition \ref{random_test} and Remark \ref{rk_random_test}. 

\begin{proposition} \label{psh_test}
Let $\varphi$ 
be a continuous p.s.h. function on $D$.
There exists a constant $M_\varphi$ such that
if $(R_m)\subset \Cc^1_v(D)$ and
$(S_n)\subset \Cc^1_h(D)$, then
$$\limsup_{m,n\rightarrow\infty}\langle d^{-m-n} (f^m)^*R_m\wedge (f^n)_*S_n,\varphi\rangle
\leq M_\varphi.$$ 
If $R_n$ and $S_n$ are as in Theorem \ref{equilibrium_measure}, we have   
$$\lim_{n\rightarrow\infty}\langle d ^{-2n}(f^n)^*R_n\wedge (f^n)_*S_n,\varphi\rangle= M_\varphi.$$ 
\end{proposition}
\begin{proof}
By Proposition \ref{intersection_explicite_construction}, we can assume that $R_n$ and 
$S_n$ are smooth forms.
We can also assume that  $m\geq n$ and $n\rightarrow\infty$.
Write $d^{-m}(f^m)^*R_m=d^{-n}(f^n)^*R_{m,n}$ with
$R_{m,n}:=d^{-m+n}(f^{m-n})^*R_m$. This allows us 
to suppose that $m=n$.

Define the currents $T_n$ in $\Cc^1_v(D^2)$ by $T_n:=R_n\otimes S_n$ and
$\widetilde\varphi(z_1,w_2,z_2,w_1):=\varphi(z_1,w_1)$. Then
$$\langle (f^n)^*R_n\wedge (f^n)_*S_n,\varphi\rangle =\langle F^{n*}(T_n), 
\widetilde\varphi[\Delta]
\rangle.$$
The current $\Phi:=\widetilde\varphi[\Delta]$ is not horizontal, but $F_*(\Phi)$ is horizontal.
Hence, Proposition \ref{psh_test} is a consequence of 
Proposition \ref{random_test} and Remark \ref{rk_random_test} applied
to $F$.
\end{proof}

We can now define the positive measure $\mu$ by
$$\langle \mu,\varphi\rangle :=M_\varphi.$$
Consider smooth forms $R\in\Cc^1_v(D')$ with support in $D'\setminus\Kc_+$ 
and $S\in\Cc^1_h(D'')$ with support in $D''\setminus\Kc_-$.
We have $\mu=\lim d^{-2n} (f^n)^*R\wedge (f^n)_*S$. Hence, $\mu$ is supported in 
the boundary of $\Kc=\Kc_+\cap \Kc_-$. Theorem \ref{th_intersection_currents} 
shows that $\mu$ is a probability measure. 

We also have
\begin{eqnarray*}
f^*(\mu) & = & \lim_{n\rightarrow\infty} d^{-2n}f^*\big((f^n)^*R\wedge (f^n)_*S\big)= 
\lim_{n\rightarrow\infty} d^{-2n}(f^{n+1})^*R\wedge (f^{n-1})_*S \\
& = & \lim_{n\rightarrow\infty} d^{-2n+2} (f^{n-1})^*(d^{-2}f^{2*}R)\wedge (f^{n-1})_*S = \mu.
\end{eqnarray*}
Hence, $\mu$ is invariant.
\\

The following corollary gives us an extremality property of $\mu$:
\begin{corollary} \label{extremality_measure}
Let $(R_m)\subset \Cc^1_v(D)$ and $(S_n)\subset \Cc^1_h(D)$.
Let $\nu$ be a limit value of
$d^{-m-n} (f^m)^*R_m\wedge (f^n)_*S_n$ when $\min(m,n)\rightarrow\infty$. 
Then 
$$\langle \nu,\varphi \rangle
\leq \langle \mu,\varphi\rangle \quad \mbox{for } \varphi \mbox{ p.s.h. on } D.$$ 
If $\varphi$ is pluriharmonic, then 
$\langle \nu,\varphi\rangle =\langle\mu,\varphi\rangle$. 
\end{corollary}
\begin{proof}
We can assume that $\varphi$ is continuous since we can approximate it by a decreasing sequence
of continuous p.s.h. functions.
Proposition \ref{psh_test} implies that $\langle \nu,\varphi\rangle\leq 
\langle \mu,\varphi\rangle$. When $\varphi$ is pluriharmonic, this inequality holds
for $-\varphi$. Hence $\langle \nu,\varphi\rangle =\langle\mu,\varphi\rangle$.
\end{proof}

The proof of the following results are left to the reader (see Corollaries \ref{random_dominating_convergence}, 
\ref{random_varieties_convergence}, \ref{one_test_form} and Proposition \ref{random_one_test_form}).

\begin{corollary} \label{dominating_convergence_measure}
Let $R_n$, $R_n'$ in $\Cc^1_v(D)$ and  $S_n$, $S_n'$ in $\Cc^1_h(D)$ 
 and
$c>0$ such that $R_n'\leq cR_n$, $S_n'\leq cS_n$ for every $n$. 
Let $(m_i)$
and $(n_i)$ be increasing sequences of integers. If 
$$d^{-m_i-n_i}(f^{m_i})^*R_{m_i}\wedge (f^{n_i})_*S_{n_i}
\rightarrow\mu,$$ 
then  $$d^{-m_i-n_i}(f^{m_i})^*R'_{m_i}\wedge (f^{n_i})_*S_{n_i}'
\rightarrow\mu.$$
\end{corollary}

\begin{proposition} \label{one_test_function}
Let $R_m$, $S_n$, $m_i$, $n_i$ be as in Corollary \ref{dominating_convergence_measure}. 
Let $\varphi$ be a function strictly 
p.s.h. on $D$. Then, 
$$d^{-m_i-n_i}(f^{m_i})^*R_{m_i}\wedge (f^{n_i})_*S_{n_i}\rightarrow\mu$$ 
if and only if 
$$\langle d^{-m_i-n_i}(f^{m_i})^*R_{m_i}\wedge (f^{n_i})_*S_{n_i},\varphi\rangle
\rightarrow\langle \mu,\varphi\rangle.$$
\end{proposition}

\begin{corollary} \label{distribution_intersection_varieties}
Let $(n_i)$ be an increasing sequence of integers. Then, there exist a subsequence
$(m_i)$ and a pluripolar set $\Ec\subset D$ such that, for every $(a,b)\in D\setminus \Ec$, we have
$$d^{-2m_i}(f^{m_i})^*[z=a]\wedge (f^{m_i})_*[w=b]\rightarrow \mu$$
where $(z,w)$ are the coordinates 
of $\C^p\times\C^{k-p}$.
\end{corollary}

To complete the proof of Theorem \ref{equilibrium_measure}, we have only to check that $\mu$ is mixing.
That is
\begin{eqnarray}
\lim_{m\rightarrow\infty} \langle\mu, (\phi\circ f^m)(\psi\circ f^{-m}) \rangle & = &  
\langle \mu,\phi \rangle \langle\mu,\psi \rangle \label{eq_for_mixing}
\end{eqnarray}
for every functions $\phi$ and $\psi$ smooth in a neighbourhood of $\overline D$. 
Define a function $\varphi$ on $D^2$ by 
$$\varphi(z_1,w_2,z_2,w_1):=\phi(z_1,w_1)\psi(z_2,w_2).$$

\begin{lemma} \label{semi_mixing}
Assume that $\varphi$ is p.s.h. Then 
$$\limsup_{m\rightarrow\infty} \langle\mu, (\phi\circ f^m)(\psi\circ f^{-m}) \rangle \leq 
\langle \mu,\phi \rangle \langle\mu,\psi \rangle.$$
\end{lemma}
\begin{proof} 
Let $R\in\Cc^1_v(D')$ and $S\in\Cc^1_h(D'')$ be smooth forms.
Define $T:=R\otimes S$ and $T'=S\otimes R$. We have 
\begin{eqnarray*}
\langle\mu, (\phi\circ f^m)(\psi\circ f^{-m}) \rangle & = & 
\lim_{n\rightarrow\infty} \langle d^{-2n} (F^n)^*T, (\varphi\circ F^m)[\Delta]\rangle \\
& = & \lim_{n\rightarrow\infty} \langle d^{-2n} (F^m)^*\big((F^{n-m})^*T
\varphi \big), [\Delta]\rangle \\
& = & \lim_{n\rightarrow\infty} \langle d^{-2n} (F^{n-m})^*T
\varphi, (F^m)_*[\Delta]\rangle \\
& = & \lim_{n\rightarrow\infty} 
\langle d^{-2n} (F^{n-m})^*T\wedge (F^m)_*[\Delta],\varphi\rangle.
\end{eqnarray*}

Applying Proposition \ref{psh_test} to  $F$ gives
\begin{eqnarray*}
\lefteqn{\limsup_{m\rightarrow\infty} 
\langle\mu, (\phi\circ f^m)(\psi\circ f^{-m}) \rangle}\\
& \leq & 
\lim_{m\rightarrow\infty} \langle d^{-4m} (F^m)^*T \wedge (F^m)_*T',\varphi \rangle\\
& = & \lim_{m\rightarrow\infty} \langle d^{-2m} (f^m)^*R\wedge (f^m)_*S,\phi\rangle 
\langle d^{-2m} (f^m)^*R\wedge (f^m)_*S,\psi\rangle\\
& = & \langle\mu,\phi\rangle \langle \mu,\psi\rangle.
\end{eqnarray*}
\end{proof}

\noindent
{\bf End of the proof of Theorem \ref{equilibrium_measure}.}
Since $\phi$ and $\psi$ can be written as differences of smooth strictly p.s.h. 
functions, in order to prove (\ref{eq_for_mixing}), it is sufficient to consider $\phi$ and $\psi$
smooth strictly p.s.h. 
in a neighbourhood of $\overline D$. Let $A>0$ be a large constant.
Then, $(\phi(z_1,w_1)+A)(\psi(z_2,w_2)+A)$ is p.s.h. Lemma \ref{semi_mixing} implies that
$$\limsup_{m\rightarrow\infty} \langle \mu, (\phi\circ f^m+A)(\psi\circ f^{-m}+A)\rangle 
\leq \langle \mu,\phi+A\rangle \langle \mu, \psi+A \rangle.$$
Since $\mu$ is invariant, we have $\langle \mu, \phi\circ f^m \rangle = 
\langle \mu,\phi \rangle$ and $\langle \mu, \psi\circ f^{-m} \rangle = 
\langle \mu,\psi \rangle$. We deduce from the last inequality that 
\begin{eqnarray}
\limsup_{m\rightarrow\infty} \langle \mu, (\phi\circ f^m)(\psi\circ f^{-m})\rangle & \leq & \langle \mu, \phi
\rangle \langle \mu,\psi\rangle. \label{ineq_for_mixing}
\end{eqnarray}

On the other hand, the function $(\phi(z_1,w_1)-A)(-\psi(z_2,w_2)+A)$ is also 
p.s.h. in a neighbourhood of $\overline D$. In the same way, we obtain
\begin{eqnarray}
\limsup_{m\rightarrow\infty} -\langle \mu, (\phi\circ f^m)(\psi\circ f^{-m})\rangle & \leq & 
-\langle \mu, \phi \rangle \langle \mu,\psi\rangle. \label{ineq_for_mixing_bis}
\end{eqnarray}
The inequalities (\ref{ineq_for_mixing}) and (\ref{ineq_for_mixing_bis}) 
imply (\ref{eq_for_mixing}). Hence, $\mu$ is mixing.
\hfill $\square$
\\

The following proposition generalizes the convergence in Theorem \ref{equilibrium_measure}.

\begin{proposition} \label{stronger_convergence_measure}
Let $(R_m)\subset \Cc^1_v(D')$ and $(S_n)\subset \Cc^1_h(D'')$ be uniformly bounded 
sequences of continuous forms. Then, $d^{-m-n}(f^m)^*R_m\wedge (f^n)_*S_n$ converges weakly 
to $\mu$ when $\min(m,n)\rightarrow\infty$.
\end{proposition}
\begin{proof} It is sufficient to consider the case where 
$m\leq n$ and $m\rightarrow\infty$. 
If we replace $M'$, $N''$ by bigger domains, we can assume that there exist $c>0$ and
continuous forms $R\in \Cc^1_v(D')$ and $S\in \Cc^1_h(D'')$ such that $R_m\leq cR$
and $S_n\leq cS$ for every $m$ and $n$. 
By Corollary \ref{dominating_convergence_measure}, it is sufficient to prove that   
$d^{-m-n}(f^m)^*R\wedge (f^n)_*S\rightarrow\mu$. 
We will use the same idea as in Theorem \ref{random_green_current}.

Let $\varphi$ be a continuous function strictly p.s.h. on $D$ with $0\leq \varphi\leq 1$. 
By Proposition \ref{one_test_function}, we only need to check that
$\langle d^{-m-n}(f^m)^*R\wedge (f^n)_*S,\varphi\rangle\rightarrow M_\varphi$. 
Write
\begin{eqnarray*}
\langle d^{-m-n}(f^m)^*R\wedge (f^n)_*S,\varphi\rangle & = & \langle R, d^{-m-n} (\varphi\circ f^{-m}) 
(f^{m+n})_*S\rangle\\
& =: & \langle R, \Psi_{m,n}\rangle.
\end{eqnarray*}
Observe that each $\Psi_{m,n}$ is positive, bounded by a current in $\Cc^1_h(D'')$ and verifies 
$\ddc\Psi_{m,n}\geq 0$. If $R_\theta$ is defined as in Section \ref{geometry_currents}, 
then $\phi_{m,n}(\theta):=
\langle R_\theta, \Psi_{m,n} \rangle$ define a uniformly bounded family of subharmonic functions on 
$\theta\in V$. Since $R=R_{1}$, we want to prove that $\phi_{m,n}(1)\rightarrow M_\varphi$.
By Proposition \ref{psh_test}
$$\limsup_{m,n\rightarrow\infty} \phi_{m,n}(\theta)\leq M_\varphi$$
and by Lemma \ref{disc2}
$$\lim_{\theta\rightarrow1}\sup_{m,n}|\phi_{m,n}(\theta)-\phi_{m,n}(1)|=0.$$ 
Hence, it is sufficient to prove that $\phi_{m,n}$ converge to $M_\varphi$ in $\Loneloc(V)$. 
By maximum principle, we only have to check that $\phi_{m,n}(0)=
\langle R_0,\Psi_{m,n}\rangle\rightarrow M_\varphi$.

Consider a smooth form $R'\in\Cc^1_v(D')$ and define $R'_{m,n}:=d^{m-n}(f^{n-m})^*R'$. Theorem 
\ref{equilibrium_measure}
implies that $d^{-m-n}(f^m)^*R'_{m,n}\wedge (f^n)_*(S)\rightarrow\mu$. Let $R'_{m,n,\theta}$ be 
the currents of the structural discs associated to $R'_{m,n}$ constructed in Section \ref{geometry_currents}. 
Then, $\phi'_{m,n}(\theta):=\langle R'_{m,n,\theta}, \Psi_{m,n}
\rangle$ define a uniformly bounded family of subharmonic functions on $\theta\in V$. We also have 
$\limsup  \phi'_{m,n}(\theta) \leq M_\varphi$ and 
$\lim  \phi'_{m,n}(1) =M_\varphi$ since $R'_{m,n,1}=R'_{m,n}$. By maximum
principle, 
$\phi'_{m,n}\rightarrow M_\varphi$ in $\Loneloc(V)$. Lemma \ref{disc1}
implies that 
$$\lim_{\theta\rightarrow 0}\sup_{m,n}|\phi'_{m,n}(\theta)-\phi_{m,n}'(0)|=0.$$
Hence, 
 $\langle R_{m,n,0}', \Psi_{m,n}\rangle=\phi'_{m,n}(0)\rightarrow M_\varphi$. 
We have seen in Section \ref{geometry_currents} that $R_0=R_{m,n,0}'$.
It follows that  $\langle R_0, \Psi_{m,n}\rangle\rightarrow M_\varphi$. 
\end{proof}

\begin{corollary} \label{convergence_on_current}
Let $S\in \Cc^1_h(D)$ be a continuous form. Then
$d^{-n}T_+\wedge (f^n)_*S$ converge weakly to $\mu$. 
\end{corollary}
\begin{proof} Let $\varphi$ be a continuous strictly p.s.h. function on $D$. 
Let $R\in\Cc^1_v(D)$ be a smooth form. Corollary \ref{green_current} 
implies that $T_+=\lim d^{-n}f^{n*}(R)$.
Hence, there exists $m>n$ such that
$$ \left|\langle d^{-m-n}(f^m)^*R\wedge (f^n)_*S,\varphi\rangle - \langle d^{-n}T_+\wedge (f^n)_*S,\varphi\rangle
\right|\leq 1/n.$$
By Proposition \ref{stronger_convergence_measure}, this implies that 
$\lim \langle d^{-n}T_+\wedge (f^n)_*S,\varphi\rangle = \langle\mu,\varphi\rangle$.
Proposition \ref{one_test_function} implies that $\lim d^{-n}T_+\wedge (f^n)_*S=\mu$.
\end{proof}

We now show that the equilibrium measure  
is equal to the wedge product of the Green currents.

\begin{theorem} \label{intersection_green_currents}
Let $f$ be an invertible  horizontal-like map and $\mu$, $T_+$, $T_-$ be 
as above. Then $$\mu=T_+\wedge T_-.$$ 
\end{theorem}
\begin{proof} 
Let $\varphi$ be a continuous p.s.h. function on $D$.
Let $R\in \Cc_v^1(D)$ and $S\in\Cc_h^1(D)$ be smooth forms. Corollary \ref{green_current} and 
Theorem \ref{equilibrium_measure} implies 
that $d^{-n}(f^n)^*R\rightarrow T_+$,  $d^{-n}(f^n)_*S\rightarrow T_-$ and 
 $d^{-2n}(f^n)^*R\wedge (f^n)_*S\rightarrow \mu$. It follows from Theorem \ref{th_intersection_currents}
that $\langle \mu,\varphi\rangle \leq \langle T_+\wedge T_-,\varphi\rangle$. 

On the other hand, we have $f^*T_+=dT_+$ and $f_* T_-= d T_-$. Hence 
Proposition \ref{psh_test} imply that 
$$\langle   T_+\wedge T_-,\varphi\rangle = \lim \langle d^{-2n} (f^n)^*T_+\wedge (f^n)_*T_-,\varphi\rangle
\leq \langle\mu,\varphi\rangle.$$
Theorem \ref{intersection_green_currents} follows.
\end{proof}

\section{Entropy}
\label{section_entropy}
We will show that the topological entropy $h_t(f_{|\Kc})$ of the 
restriction of $f$ to the invariant compact set $\Kc$ is
equal to $\log d$. From the variational principle \cite{Walters}, 
it follows that the entropy of $\mu$ is bounded from above by $\log d$.
We will show that this measure has entropy $h(\mu)=\log d$. This also implies that 
$h_t(f_{|\supp(\mu)})=\log d$.

\begin{theorem} \label{entropy}
Let $f$, $d$, $\Kc$, $\mu$ be as above. Then, the topological entropy of $f_{|\Kc}$ is equal to
$\log d$ and $\mu$ is an invariant measure of maximal entropy $\log d$.
\end{theorem}

We have to prove that $h_t(f_{|\Kc})\leq\log d$ and $h(\mu)\geq \log d$. 
Using Yomdin's results \cite{Yomdin},
Bedford-Smillie proved the second inequality for H{\'e}non maps
\cite{BedfordSmillie} (see also 
Smillie \cite{Smillie}). 
We only need the following lemma applied to a closed form
$S$ strictly positive in a neighbourhood 
of $\Kc_-$
in order to adapt their proof and get $h(\mu)\geq \log d$.

\begin{lemma} \label{lemma_for_entropy}
Let $S\in \Cc^1_h(D'')$ be a smooth 
form. 
Then, there exist an increasing sequence $(n_i)$ of positive integers and a point $a\in M'$
such that 
$$\frac{1}{n_i}\sum_{j=0}^{n_i-1} d^{-n_i} (f^j)^*[z=a]\wedge (f^{n_i-j})_*S\rightarrow\mu.$$ 
\end{lemma}
\begin{proof} Let $\varphi$ be a smooth strictly p.s.h. function on $D$. 
Define a sequence of p.s.h. functions, for $a\in M$ (see Theorem \ref{th_slices}):
$$\phi_n(a):=\left\langle\frac{1}{n}\sum_{j=0}^{n-1} d^{-n} (f^j)^*[z=a]\wedge (f^{n-j})_*S, 
\varphi \right\rangle.$$ 

Let $\nu$ be a smooth probability measure on $M'$. Consider the smooth form 
$R:=\pi_1^*(\nu)$ in $\Cc^1_v(D')$. Proposition \ref{stronger_convergence_measure} implies that 
$$\frac{1}{n}\sum_{j=0}^{n-1} d^{-n} (f^j)^*R\wedge (f^{n-j})_*S\rightarrow\mu.$$ 
Hence, $\int\phi_n(a)\d\nu(a)\rightarrow M_\varphi$. 
On the other hand, Proposition \ref{psh_test} implies that 
$$\lim_{n\rightarrow\infty}\sup_{a\in M'} \phi_n(a)\leq M_\varphi.$$
It follows that there exist $(n_i)$ and $a\in M'$
such that $\lim \phi_{n_i}(a)= M_\varphi$. 
As in Propositions \ref{random_one_test_form} and \ref{one_test_function}, 
we prove that $(n_i)$ and $a$ satisfy  the lemma.
\end{proof}

Now, we prove the first inequality $h_t(f_{|\Kc})\leq\log d$. Analogous inequalities have been proved in 
\cite{Gromov, DinhSibony1, DinhSibony2, DinhSibony3}. We use here some arguments in Gromov \cite{Gromov}
and in \cite{DinhSibony1}. 

Let $\Gamma_{[n]}$ be the graph of the map $x\mapsto (f(x),\ldots, f^{n-1}(x))$. 
This is the set of points $(x,f(x),\ldots,f^{n-1}(x))$.
We use the canonical euclidian metric on $D^n$. 
Let $D_*:=M'\times N''$. We have $\Kc\subset D_*\Subset D$. Define 
$$\lov(f):=\limsup_{n\rightarrow\infty} \frac{1}{n}\log\vol(\Gamma_{[n]}\cap D_*^n).$$
Following Gromov \cite{Gromov, DinhSibony1}, 
we have $h_t(f_{|\Kc})\leq \lov(f)$. We will show that $\lov(f)\leq\log d$; 
then  $h_t(f_{|\Kc})=\lov(f)=\log d$ since $h_t(f_{|\Kc})\geq h(\mu)\geq\log d$. 

Let $\Pi$ denote the projection of $D^n=(M\times N)^n$ on the product $M\times N$ of the last factor $M$ and 
the first factor $N$.
Let $\Pi_1$ (resp. $\Pi_2$) denote the projections of $D^n$ 
on the product $M^{n-1}$ (resp. $N^{n-1}$) 
of the other factors $M$ (resp. $N$). Observe that
$\Pi:\Gamma_{[n]}\rightarrow M\times N$ is proper 
and defines a ramified covering of degree $d^{n-1}$ over $M\times N$.
Indeed, for a generic point $(a,b)\in M\times N$ the fiber
$\Pi^{-1}(a,b)\cap\Gamma_{[n]}$ contains a number of points equal to
the number of points in $\{z=a\}\cap f^{n-1}\{w=b\}$, i.e. equal to
$d^{n-1}$ (see Proposition \ref{operators_on_currents}).
Moreover, we have 
$\Gamma_{[n]}\subset \Pi_1^{-1}({M'}^{n-1})$ and $\Gamma_{[n]}\subset \Pi_2^{-1}({N''}^{n-1})$. Now, 
it is sufficient
to apply the following lemma (see \cite[lemme 3.3.3]{DinhSibony1} for the proof).

\begin{lemma} \label{volume_estimate}
Let $\Gamma$ be an analytic subset of dimension $k$ of $D\times M^m\times N^m$
such that $\Gamma\subset D\times {M'}^m\times {N''}^m$. 
We assume that $\Gamma$ is a ramified covering over $D$ of degree 
$d_\Gamma$.
Then, there exist $c>0$, $s>0$ independent of $\Gamma$ and of $m$ such that
$$\vol(\Gamma\cap D_*\times M^m\times N^m)\leq cm^sd_\Gamma.$$ 
\end{lemma}

Tien-Cuong Dinh \hfill Nessim Sibony\\
Institut de Math{\'e}matique de Jussieu \hfill  Math{\'e}matique - B{\^a}timent 425 \\
Plateau 7D, Analyse Complexe \hfill  UMR 8628\\
175 rue du Chevaleret \hfill  Universit{\'e} Paris-Sud \\
75013 Paris, France \hfill 91405 Orsay, France \\
{\tt dinh@math.jussieu.fr} \hfill {\tt nessim.sibony@math.u-psud.fr} \\

\end{document}